\RequirePackage{fix-cm}
\documentclass[smallextended]{svjour3}       
\smartqed  
\usepackage{graphicx}
\usepackage{amsmath}
\usepackage{amsfonts}
\usepackage{bm}
\usepackage{bbm}
\usepackage{comment}
\usepackage{hyperref}
\usepackage{soul,color}
\usepackage[titletoc]{appendix}
\usepackage[square,sort&compress,comma,numbers]{natbib} 
\usepackage{amsmath,amsfonts,amssymb,theorem}
\usepackage{url,subfig,anysize,wrapfig}
\usepackage{multirow,tabularx,afterpage,fancyhdr}
\usepackage{cases}
\usepackage{epstopdf}
\usepackage{bm}

\usepackage[linewidth=1pt]{mdframed}
\usepackage{lipsum}

\usepackage{lineno}
\usepackage{placeins}
\usepackage[pointedenum]{paralist}
\usepackage{mdwlist}

\usepackage{tikz}
\usetikzlibrary{calc}

\usepackage{todonotes}
\usepackage{xcolor}

\newcommand{\mydrafttext}{}
\newcommand{\drafttext}[1]{\renewcommand{\mydrafttext}{#1}}
\newboolean{draft}  
\setboolean{draft}{true}
\ifthenelse{\boolean{draft}}
{
    \newcounter{comments}
    \drafttext{{\color{red}This document is a draft.}}
    \newcommand{\shane}[1]{\addtocounter{comments}{1}{\color{red}[Shane comment \thecomments: #1]}}
    \newcommand{\mattia}[1]{\addtocounter{comments}{1}{\color{blue}[Mattia comment \thecomments: #1]}}
    \newcommand{\peter}[1]{\addtocounter{comments}{1}{\color{orange}[Peter comment \thecomments: #1]}}
}
{
\newcommand{\shane}[1]{}
\newcommand{\mattia}[1]{}
\newcommand{\peter}[1]{}
}

\journalname{}
\begin{document}

\title{Global invariant manifolds delineating transition and escape dynamics in dissipative systems}

\titlerunning{Global invariant manifolds delineating transition and escape dynamics in dissipative systems}
\authorrunning{Zhong and Ross} 

\author{Jun Zhong \and Shane D. Ross
}

\institute{Jun Zhong (Corresponding author) \at
	School of Engineering, Brown University,
	Providence, RI 02912, USA\\
	\email{jun\_zhong@brown.edu}    \\      
	\\
	Shane D. Ross \at
	Aerospace and Ocean Engineering, Virginia Tech,
	Blacksburg, VA 24061, USA \\
	\email{sdross@vt.edu}
}

\date{This version: \today}


\maketitle

\begin{abstract}
Invariant manifolds play an important role in organizing global dynamical behaviors. For example, it is found that in multi-well conservative systems where the potential energy wells are  connected by index-1 saddles, the motion between potential wells is governed by the invariant manifolds of a periodic orbit around the saddle. In two degree of freedom systems, such invariant manifolds appear as cylindrical conduits which are referred to as transition tubes. In this study, we apply the concept of invariant manifolds to study the transition between potential wells in not only conservative systems, but more realistic  dissipative systems, by solving respective proper boundary-value problems. The example system considered is a two mode model of the snap-through buckling of a shallow arch. We define the transition region, $\mathcal{T}_h$, as a set of initial conditions of a given initial Hamiltonian energy $h$ with which the trajectories can escape from one potential well to another, which in the example system corresponds to  snap-through buckling of a structure. The numerical results reveal that in the conservative system the boundary of the transition region, $\partial \mathcal{T}_h$, is a cylinder, while in the dissipative system, $\partial \mathcal{T}_h$ is an ellipsoid. The algorithms developed in the current research from the perspective of invariant manifold provides a robust theoretical-computational framework to study  escape and transition dynamics.


\end{abstract}

\section{Introduction}
Escape or transition between potential wells is  found in a  number of important systems, such as snap-through buckling of curved structures  \cite{zhong2018tube,collins2012isomerization}, chemical reactions \cite{OzDeMeMa1990,DeMeTo1991,wiggins2001impenetrable,uzer2002geometry,gabern2005theory}, celestial mechanics \cite{jaffe2002statistical,KoMaRoLoSc2004,onozaki2017tube},
and capsize of floating structures \cite{sequeira2018manifestation,NaRo2017}, to name but a few. The prediction of transition plays an important role in both utilization and evasion. In a one degree of freedom system, the only possible escape route out of a well is via a local hilltop of the potential energy. The situation becomes more complicated for higher degree of freedom systems, since there are infinitely many routes to escape, generally via index-1 saddles connecting the potential wells. 
The transition boundary for the possible escape trajectories in conservative higher dimensional systems has been demonstrated to be the stable invariant manifold of a normally hyperbolic invariant manifold (NHIM)  around the index-1 saddle \cite{uzer2002geometry}. 
For a given energy, the  trajectories inside the stable manifold 
can escape from the potential well, while those just outside the stable manifold bounce back from the saddle and do not escape the potential well. Thus, the general way of computing the transition boundary in conservative systems is to find the NHIM at a fixed energy associated to the index-1 saddle and then compute its stable and unstable invariant manifolds. In the two degree of freedom case, the NHIM is a collection of periodic orbits, where each orbit corresponds to a choice of energy, and the corresponding stable manifold is geometrically cylindrical and sometimes called a transition tube.

\paragraph{Mechanical systems with dissipation.}
However, the situation for a system with  dissipation added is  less well-understood. When energy dissipation is incorporated, the bound orbits comprising a NHIM at constant energy no longer exist which makes the classical method of computing the transition boundary fail for dissipative systems. 
Thus, a new framework must be established for  dissipative systems. Ref.\ \cite{zhong2018tube} proposed a bisection method to find the transition boundary on a specific Poincar\'e section for both conservative and dissipative systems. The whole phase space structure that governs the transition was not discussed, although it can be obtained by selecting a collection of Poincar\'e sections. This method is versatile, but 
can be inefficient 
if too many Poincar\'e sections are considered or the shape of the transition boundary on a Poincar\'e section is irregular and distorted. 
The authors \cite{zhong2020geometry} summarized the geometry of phase space structures that governs the escape from potential wells in some widely known systems with two degrees of freedom where gyroscopic and dissipative forces have been added. 
They found that the transition boundary for  a specific given initial energy goes from a {\it cylindrical} tube in the conservative system to an {\it ellipsoid} in the dissipative system, referred to as the transition tube and transition ellipsoid, respectively. 
While the transition tube is the stable invariant manifold of the periodic orbit of the initial energy in the conservative system, the transition ellipsoid is a subset of the stable invariant manifold of  index-1 saddle equilibrium point in the dissipative system. 
The previous paper discussed the linearized dynamics \cite{zhong2020geometry}. 
These topological results carry over to the nonlinear setting via the stable manifold theorem \cite{meiss2007differential,wiggins2003introduction,perko2013differential} and a theorem of Moser \cite{Moser1958,Moser1973}, for the dissipative and conservative cases, respectively.
This paper aims to extend to the nonlinear setting using the invariant manifold perspective.

\paragraph{Invariant manifolds.}
The concept of invariant manifolds is crucial for understanding the characteristics of a dynamical system. In general, the global invariant manifold cannot be obtained analytically. Thus, numerical and computational algorithms are necessary. Suppose we have a continuous dynamical system written as a set of autonomous ordinary differential equations,
\begin{equation}
\dot{x}=f(x),
\label{general ODEs}
\end{equation}
where $x\in \mathbb{R}^n$ and the vector field $f: \mathbb{R}^n \rightarrow \mathbb{R}^n$ is sufficiently smooth. Here the dot over the quantity is the derivative with respect to time $t\in \mathbb{R}$. For the dissipative mechanical systems envisioned, the vector field is assumed to have a hyperbolic equilibrium point at $x=x_e$, i.e.,  $f(x_e)=0$. Its Jacobian matrix $Df(x_e)$ has $k$ eigenvalues with negative real part and $n-k$ eigenvalues with positive real part. The eigenvectors corresponding to the eigenvalues with negative and positive real parts are denoted by $u_i$ and $v_i$, respectively. Thus, the spaces spanned by $u_i$ and $v_i$ are referred to stable and unstable subspaces of the linearized system, denoted by $E^s$ and $E^u$, which are defined by,
\begin{equation}
\begin{aligned}
E^s&= \text{span} \{u_1,u_2, \cdots, u_k\},\\
E^u&= \text{span} \{v_1,v_2,\cdots, v_{n-k}\}.
\end{aligned}
\end{equation}
From the Theorem of the Local Stable and Unstable Manifold \cite{meiss2007differential,wiggins2003introduction,perko2013differential}, there exists a $k$-dimensional invariant local stable manifold and a $(n-k)$-dimensional invariant unstable manifold, denoted by $W^s_{loc}(x_e)$ and $W^u_{loc}(x_e)$, which are tangent to $E^s$ and $E^u$ at $x_e$, respectively. Thinking in terms of computation, after the local stable and unstable invariant manifolds are established, the global stable and unstable invariant manifolds can be grown from the corresponding local invariant manifold \cite{meiss2007differential,wiggins2003introduction,perko2013differential,krauskopf2006survey} which are defined by,
\begin{equation}
\begin{aligned}
W^s(x_e)&=\left\{x\in \mathbb{R}^n ~| \lim_{t \rightarrow +\infty} \phi_t(x)=x_e \right\}= \bigcup_{t \geqslant 0} \phi_t(W^s_{loc}(x_e)),\\
W^u(x_e)&=\left\{x\in \mathbb{R}^n ~| \lim_{t \rightarrow -\infty} \phi_t(x)=x_e \right\} = \bigcup_{t \leqslant 0} \phi_t(W^u_{loc}(x_e)),
\label{global manifold}
\end{aligned}
\end{equation} 
where $\phi_t$ is the flow map of the   system \eqref{general ODEs}. From the definitions of the invariant manifolds in \eqref{global manifold}, it is intuitive to compute the global invariant manifold by numerical integration using a collection of initial conditions on a $(k-1)$-dimensional hyper-sphere with a small radius $\delta$ centered at $x_e$ in the corresponding subspace. This idea works well for a one-dimensional invariant manifold of the equilibrium point embedded in any dimensional space \cite{parker2012practical}. However, some challenges \cite{krauskopf2003computing} may appear when computing  higher dimensional invariant manifolds, such as  large aspect ratios of the manifold surface due to the significant differences in the magnitude of the real part of the eigenvalues leading trajectories on the manifold to be attracted to the most stable direction \cite{osinga2018understanding}. In this case, directly growing the invariant manifold from the local initial sphere is impractical. Other methods can be found in a review paper \cite{krauskopf2006survey} about computing the global invariant methods and interested readers are referred therein for more details.

\paragraph{Computing global invariant manifolds.}
Ref. \cite{krauskopf2003computing} presents the approach  of
computing the global invariant manifolds of a hyperbolic point as a family of orbit segments,  solving a suitable two-point boundary-value problem (BVP).
It was applied to some examples to compute a two-dimensional invariant manifold formed by a family of geodesic level sets. 
Starting from another perspective, similar to cell-mapping method, a  box covering technique \cite{dellnitz1997subdivision,dellnitz1996computation} was developed to compute  invariant manifolds. 
In this approach, a subdivision algorithm is used to produce the local invariant manifold first and then a box-oriented continuation technique is applied to extend it to the global manifold. 
Theoretically this technique is applicable to compute  invariant manifolds of arbitrary dimension.
However, due to the large number of boxes in high-dimensional systems which slows down the computation, only moderate dimensional problems are considered in practice. 
The Lagrangian descriptor \cite{madrid2009distinguished,mendoza2010hidden} is a trajectory-based diagnostic method, originally developed in the context of transport in fluid mechanics, to detect  invariant manifolds and invariant manifold-like structures \cite{naik2019finding,mancho2013lagrangian}. 
It measures the geometrical properties of particle trajectories, such as the arc-length, within a fixed forward and backward time starting at given initial conditions. Since it is an integration method, its computational expenses still need to be examined.

\paragraph{Transition region boundary across an index-1 saddle with dissipation computed via its stable global invariant manifold.}
In this study, we present an approach that systematically identifies all the trajectories which will cross from one side of an index-1 saddle to the other, for instance from one potential well to another. 
While we do this in the context of a specific physical situation  (the snap-through buckling of a shallow arch), the approach is quite general. 
We define the transition region in phase space, $\mathcal{T}_h$, as the set of initial conditions of a given initial Hamiltonian energy $h$ with which the trajectories can escape from one potential well to another.
The numerical results reveal that in the two degree of freedom conservative version of the system, the boundary of the transition region, $\partial \mathcal{T}_h$, is topologically a cylinder, while in the dissipative system, $\partial \mathcal{T}_h$ is topologically a sphere, that is, a cylinder with the `ends closed'.
The one-parameter family of boundaries,  $\partial \mathcal{T}_h$ with parameter $h\ge 0$,  can be obtained by solving a proper BVP. 
The BVP is implemented using the numerical package COCO \cite{dankowicz2013recipes}.

\section{Example two degree of freedom system: Snap-through of a shallow arch}
We use a discretized model of the snap-through buckling of a continuous shallow arch to illustrate the computation the transition boundary from the perspective of a stable invariant manifold. 
As will be seen below, the resulting potential energy surface is topologically equivalent to several other systems, including chemical systems. 

We consider an arch of length $L$, width $b$, and thickness $d$. For a slender arch, it allows us to adopt the Euler-Bernoulli beam theory taking account the von K\'arm\'an-type geometrical nonlinearity \cite{zhong2016analysis} to characterize the nonlinear deformations. Considering in-plane immovable constraints at each end of the arch, the deflection of the arch is governed by an integro-differential equation \cite{WiVi2016,zhong2018tube},
\begin{equation}
\begin{aligned}
\rho A \frac{\partial^2 w}{\partial t^2} + c_d \frac{\partial w}{\partial t} + EI \left(\frac{\partial^4 w}{\partial x^4} - \frac{\partial^4 w_0}{\partial x^4}\right) + \left[N_T - \frac{EA}{2L}  \int_0^L \left( \left(\frac{\partial w}{\partial x} \right)^2 - \left(\frac{\partial w_0}{\partial x} \right)^2 \right) \mathrm{d} x \right] \frac{\partial^2 w}{\partial x^2} =0,
\label{PDEs}
\end{aligned}
\end{equation}
where $w$ and $w_0$ are the transverse displacement and initial deflection (or imperfection) of the arch, respectively; $A$ and $I$ are the cross-sectional area and moment of inertia which result in axial stiffness and bending stiffness, denoted by $EA$ and $EI$, respectively; $\rho$ is the mass density. Due to the immovable ends, the external axial force, which can be introduced to control the initial equilibrium shapes, can not be applied. It will be replaced by the thermal loading in this study, denoted by $N_T$. Finally, $c_d$ is the coefficient of the linear viscous damping. In this analysis, we consider a clamped-clamped arch whose boundary conditions are given by.
\begin{equation}
w=0 \quad \text{and} \quad  \frac{\partial w}{\partial x}=0 \quad \text{at} \quad x=0, L
\end{equation}

Typically, the symmetric snap-through of a shallow arch is a classic example of saddle-node bifurcation in a slender arch, while in the case of an arch which is not shallow, the snap-through is asymmetric corresponding to a subcritical pitchfork bifurcation \cite{zhong2018tube,zhong2020differential,virgin2017geometric,WiVi2016,harvey2015coexisting}. To capture the asymmetrical deformation, a two-mode truncation is utilized. Thus, the deformation and the initial imperfection can be written by,
\begin{equation}
\begin{aligned}
w(x,t)&= X(t) \phi_1(x) +Y(t) \phi_2(x),\\
w_0(x)&= \gamma_1 \phi_1(x) + \gamma_2 \phi_2(x),
\end{aligned}
\label{2_mode_shapes}
\end{equation}
where $\phi_1$ and $\phi_2$ are the first two mode shapes with the following forms,
\begin{equation}
\begin{aligned}
&\phi_n = \alpha_n \left[\sinh  \frac{\beta_n x}{L} - \sin  \frac{\beta_n x}{L} + \delta_n \left(\cosh  \frac{\beta_n x}{L} - \cos  \frac{\beta_n x}{L} \right) \right],\\
&\delta_n = \frac{\sinh \beta_n - \sin \beta_n}{\cos \beta_n - \cosh \beta_n},\\
&\cos \beta_n \cosh \beta_n =1,\\
& \alpha_1 = -0.6186, \ \ \  \alpha_2 = -0.6631.
\label{mode:Virgin}
\end{aligned}
\end{equation}

Substituting the two-mode shape approximation in \eqref{2_mode_shapes} to \eqref{PDEs}, and applying the Galerkin method, one can obtain the following nonlinear ordinary differential equations,
\begin{equation}
\begin{aligned}
& M_1 \ddot X + C_1 \dot X + K_1 \left(X - \gamma_1 \right) - N_T G_1 X - \frac{EA}{2L}G_1^2 \left(\gamma_1^2 X -X^3 \right) - \frac{EA}{2L} G_1 G_2 \left(\gamma_2^2 X -X Y^2 \right)=0,\\
& M_2 \ddot Y + C_2 \dot Y + K_2 \left(Y - \gamma_2 \right) - N_T G_2 Y - \frac{EA}{2L}G_2^2 \left(\gamma_2^2 Y -Y^3 \right) - \frac{EA}{2L} G_1 G_2 \left(\gamma_1^2 Y -X^2 Y \right)=0,
\label{odes}
\end{aligned}
\end{equation}
where the coefficients are defined by,
\begin{equation}
\begin{aligned}
\left(M_i, C_i \right) = \left(\rho A, c_d \right) \int_0^L \phi_i^2 \mathrm{d}x, \ \ K_i = EI \int_0^L \left(\frac{\partial^2 \phi_i}{\partial x^2} \right)^2 \mathrm{d}x, \ \ G_i = \int_0^L \left(\frac{\partial \phi_i}{\partial x} \right)^2 \mathrm{d}x.
\label{Galerkin-coefficient}
\end{aligned}
\end{equation}

Note that \eqref{odes} can also be obtained from Lagrange's equations,
\begin{equation}
\frac{d}{dt}\left( \frac{\partial \mathcal{L}}{ \partial \dot q_i}\right) -\frac{\partial \mathcal{L}}{ \partial  q_i} = - C_i \dot q_i, \quad i =1, 2,
\end{equation}
where $q_1 = X$ and $q_2=Y$, and the Lagrangian function is given by,
\begin{equation}
\mathcal{L}(X,Y,\dot X,\dot Y) = \mathcal{T}(\dot X, \dot Y) - \mathcal{V}(X, Y),
\end{equation}
where $\mathcal{T}$ and $\mathcal{V}$ represent the kinetic energy and potential energy, respectively, given by,
\begin{equation}
\begin{split}
\mathcal{T}(\dot X, \dot Y)=& \frac{1}{2} M_1 \dot X^2 + \frac{1}{2} M_2 \dot Y^2,\\
\mathcal{V}(X, Y)=& - K_1 \gamma_1 X - K_2 \gamma_2 Y + \frac{1}{2} K_1 X^2 + \frac{1}{2} K_2 Y^2 - \frac{1}{2} N_T\left( G_1 X^2 +  G_2 Y^2 \right) \\
& - \frac{EA}{2L}G_1^2 \left(\frac{1}{2}\gamma_1^2 X^2 - \frac{1}{4}X^4 \right) - \frac{EA}{2L}G_2^2 \left(\frac{1}{2} \gamma_2^2 Y^2 - \frac{1}{4}Y^4 \right) \\
& - \frac{EA}{2L} \frac{G_1 G_2}{2} \left(\gamma_2^2 X^2 + \gamma_1^2 Y^2 -X^2 Y^2 \right).\\
\end{split}
\end{equation}

\begin{figure}[!h]
	\begin{center}
		\includegraphics[width=\textwidth]{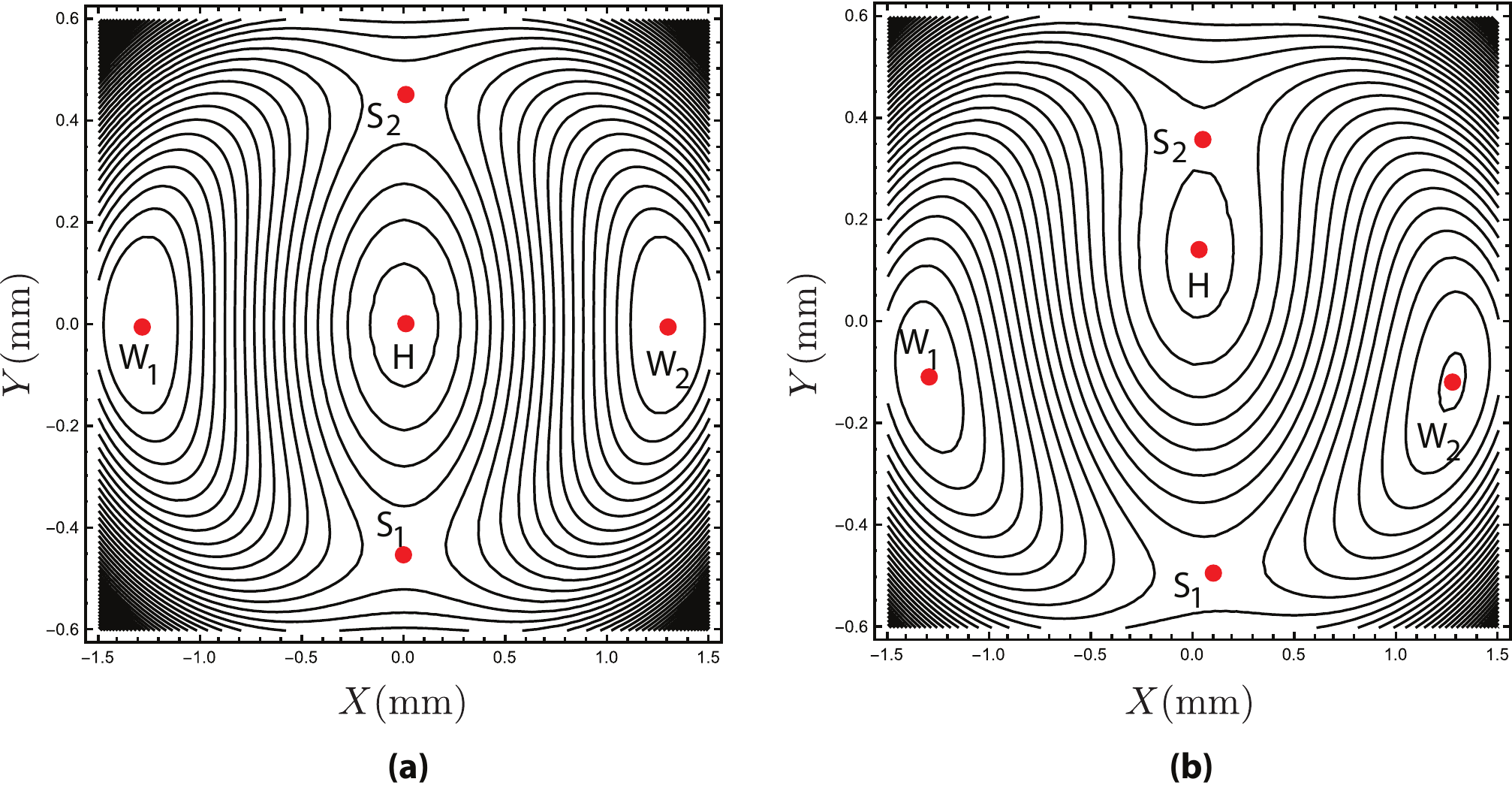}
	\end{center}
	\caption{\footnotesize
		Contours of potential energy: (a) the symmetric system, $\gamma_1=\gamma_2=0$, (b) with small initial imperfections in both modes, i.e., $\gamma_1$ and $\gamma_2$ are nonzero. Comparing the two cases, we notice the introduction of the initial imperfections changes the contours of potential energy from symmetrical to asymmetrical. 
	}
	\label{energy surface}
\end{figure}

\paragraph{Potential energy surface and links with chemistry.}
A typical form of the potential energy surface is shown Figure \ref{energy surface}(a), where W$_1$ and W$_2$ are within the two stable wells; S$_1$ and S$_2$ two index-1 saddle points; H the unstable hilltop (an index-2 saddle). 
For an equilibrium state, the system might be at rest in a position of stable equilibrium, such as point W$_1$. If the system is given a perturbation, for example, an  impact force,  snap-through buckling might occur and it can transition to the the remote equilibrium at point W$_2$. 
In general, the motion between the potential wells most likely to occur via the low-energy routes via  saddle S$_1$ or S$_2$, typically avoiding H, since H is a  potential energy maximum. 
When a small geometrical imperfection in both modes is incorporated, the symmetry of the potential energy surface is broken, as shown in Figure \ref{energy surface}(b). In the numerical examples, we will consider an imperfect shallow arch.

The potential energy surface shown is topologically equivalent to that in several two degree of freedom problems in chemistry, namely 
isomerization \cite{murrell1968symmetries,DeLi1994,collins2012isomerization}, 
double proton transfer \cite{delavega1982role,minyaev1994reaction,smedarchina2007correlated,accardi2010synchronous}, 
and other chemical reactions \cite{ezra2009phase}.

\paragraph{Hamiltonian formulation with dissipation.}
Instead of using the equations of motion from a Lagrangian perspective, in the following analysis, we put the problem to a Hamiltonian system which automatically gives first-order ordinary differential equations. Thus, we define the generalized momenta,
\begin{equation}
\begin{split}
p_i = \frac{\partial \mathcal{L}}{ \partial \dot q_i} = M_i \dot q_i,
\end{split}
\end{equation}
so $p_X = M_1 \dot X$ and $p_Y = M_2 \dot Y$,
in which case, the kinetic energy is
\begin{equation}
\mathcal{T}(p_X,p_Y)  = \frac{1}{2 M_1} p_X^2 +  \frac{1}{2 M_2} p_Y^2 ,
\end{equation}
and the Hamiltonian is 
\begin{equation}
\mathcal{H}(X,Y,p_X,p_Y) = \mathcal{T}(p_X,p_Y) + \mathcal{V}(X, Y),
\end{equation}
and Hamilton's equations (with damping) \cite{Greenwood2003} are
\begin{equation}
\begin{split}
\dot X &= \frac{\partial \mathcal{H}}{\partial p_X}=\frac{p_X}{M_1},  \hspace{0.4in} \dot Y = \frac{\partial \mathcal{H}}{\partial p_Y}=\frac{p_y}{M_2}, \\
\dot p_X &= - \frac{\partial \mathcal{H}}{\partial X} - C_H p_X=- \frac{\partial \mathcal{V}}{\partial X} - C_H p_X, \\
\dot p_Y &= - \frac{\partial \mathcal{H}}{\partial Y} - C_H p_Y=- \frac{\partial \mathcal{V}}{\partial Y} - C_H p_Y, 
\label{eq:eomHam}
\end{split}
\end{equation}
where
\begin{equation}
\begin{split}
\frac{\partial \mathcal{V}}{\partial X}=& K_1 \left(X - \gamma_1 \right) - N_T G_1 X - \frac{EA}{2L}G_1^2 \left(\gamma_1^2 X -X^3 \right) - \frac{EA}{2L} G_1 G_2 \left(\gamma_2^2 X -X Y^2 \right),\\
\frac{\partial \mathcal{V}}{\partial Y}=& K_2 \left(Y - \gamma_2 \right) - N_T G_2 Y - \frac{EA}{2L}G_2^2 \left(\gamma_2^2 Y -Y^3 \right) - \frac{EA}{2L} G_1 G_2 \left(\gamma_1^2 Y -X^2 Y \right),
\end{split}
\end{equation}
and $C_H=C_1/M_1=C_2/ M_2$ is the damping coefficient in the Hamiltonian system which can be easily found by comparing \eqref{odes} and \eqref{eq:eomHam}, and using the relations of $M_i$ and $C_i$ in  \eqref{Galerkin-coefficient}.

\section{Linearized dynamics around the equilibrium region}

\paragraph{Linearization near the index-1 saddle.}
As mentioned before, the transition between the two potential wells usually occurs around the index-1 saddles. The comprehension of local behaviors around such equilibria is an essential step to understand the transition in the complicated nonlinear system. For the geometrical and material parameters used later in the numerical examples, S$_1$ has lower potential energy than S$_2$. Thus, the potential energy of S$_1$ determines the critical energy that allows the existence of transition between the two wells. Here we focus on analyzing the linearized dynamics of S$_1$. We denote the position of S$_1$ by $x_e=(X_e,Y_e,0,0)^T$ and the linearized equations about S$_1$ can be given by,
\begin{equation}
\begin{aligned}
\dot x&= \frac{p_x}{M_1}, \hspace{0.4in } \dot y= \frac{p_y}{M_2},\\
\dot p_x&= A_{31} x + A_{32} y - C_H p_x,\\
\dot p_y&= A_{32} x + A_{42} y - C_H p_y,
\label{linearization}
\end{aligned}
\end{equation}
where $(x,y,p_x,p_y)^T= (X,Y,p_X,p_Y)^T - x_e$ and,
\begin{equation}
\begin{split}
& A_{31}= -K_1 + N_T G_1 + \frac{E A G_1^2 \left(\gamma_1^2 - 3 X_e^2 \right)}{2L} + \frac{E A G_1 G_2 \left(\gamma_2^2 -Y_e^2 \right)}{2L},\\
& A_{32}= - \frac{E A G_1 G_2 X_e Y_e}{L},\\
& A_{42}= -K_2 + N_T G_2 + \frac{E A G_2^2 \left(\gamma_2^2 - 3 Y_e^2 \right)}{2L} + \frac{E A G_1 G_2 \left(\gamma_1^2 -X_e^2 \right)}{2L}.
\label{lin paras}
\end{split}
\end{equation}

We introduce the following non-dimensional quantities,
\begin{equation}
\begin{split}
&\left( L_x, L_y \right)= L  \left(1,\sqrt{ \frac{M_1}{M_2}} \right), \omega_0= \frac{ \sqrt{- A_{32} }}{ \left( M_1 M_2\right)^ \frac{1}{4}}, \tau= \omega_0 t,  \left(\bar q_1 , \bar q_2 \right)= \left( \frac{x}{L_x}, \frac{y}{L_y} \right),\\
&\left(\bar p_1 , \bar p_2 \right)= \frac{1}{\omega_0} \left(  \frac{p_x}{ L_x M_1},\frac{p_y}{ L_y M_2}\right), \left( c_x  , c_y \right)= \frac{1}{ \omega_{0}^2} \left( \frac{A_{31}}{M_1}, \frac{A_{42}}{M_2} \right), c_1= \frac{C_H}{\omega_0}.
\label{dimless quan}
\end{split}
\end{equation} 
Making use of the non-dimensional quantities in \eqref{dimless quan}, we can rewrite the linearized equations in a non-dimensional form,
\begin{equation}
\begin{aligned}
\dot {\bar q}_1 &= \bar p_1, \hspace{0.4in} \dot {\bar q}_2 = \bar p_2,\\
\dot {\bar p}_1 &= c_x \bar q_1 - \bar q_2 - c_1 \bar p_1,\\
\dot {\bar p}_2 &=  -  \bar q_1 + c_y \bar q_2 - c_1 \bar p_2.
\label{nond eq}
\end{aligned}
\end{equation}
Written in matrix form, with column vector $\bar z=(\bar q_1 , \bar q_2 , \bar p_1 , \bar p_2)^T$, we have,
\[
\dot {\bar z} = A \bar z + D \bar z,
\]
where,
\begin{equation}
A = \begin{pmatrix}
0     	 & 0 & 1 & 0 \\
0	  & 0 & 0 & 1 \\
c_x	  & -1 & 0 & 0 \\
-1	  & c_y & 0 & 0 
\end{pmatrix},
\hspace{0.5in}
D = \begin{pmatrix}
0     	 & 0 & 0 & 0 \\
0	  & 0 & 0 & 0 \\
0	  & 0 & -c_1 & 0 \\
0	  & 0 & 0 & -c_1 
\end{pmatrix},
\label{A_and_D_matrix}
\end{equation}
are the Hamiltonian part and damping part of the linear equations, respectively. The resulting quadratic Hamiltonian is written by,
\begin{equation}
\mathcal{H}_2= \tfrac{1}{2}\bar p_1 ^2 + \tfrac{1}{2}\bar p_2 ^2 
- \tfrac{1}{2}c_x \bar q_1 ^2 - \tfrac{1}{2}c_y \bar q_2 ^2 + \bar q_1 \bar q_2,
\label{H_2_bar}
\end{equation}

\subsection{Conservative system}

\label{linear_conservative_system}
In this part, we discuss the linearized dynamics in the conservative system (i.e., $c_1=C_H=0$). It is straightforward to find that the eigenvalues of the linearized conservative system have the form $\pm \lambda$ and $\pm i \omega_p$, where $\lambda$ and $\omega_p$ are positive real numbers with the following forms,
\begin{equation}
\lambda=\sqrt{\alpha_1}, \hspace{0.2in} \omega_p=\sqrt{-\alpha_2}, \hspace{0.2in} \text{where } \alpha_{1,2}=\frac{c_x + c_y \pm \sqrt{(c_x - c_y)^2 + 4} }{2}.
\end{equation}
The corresponding generalized eigenvalues are, 
\begin{equation}
\begin{aligned}
u_{\omega_p}&=(1,c_x+\omega_p^2,0,0)^T, \hspace{0.2in }&& v_{\omega_p}=(0,0,\omega_p,c_x \omega_p + \omega_p^3)^T,\\
u_{+\lambda}&=(1,c_x-\lambda^2,\lambda,c_x \lambda-\lambda^3)^T, \hspace{0.2in} &&
u_{-\lambda}=-(1,c_x - \lambda^2, -\lambda,\lambda^3-c_x \lambda)^T.
\end{aligned}
\end{equation}

To better understand the dynamics in the phase space, we introduce a linear change of
coordinates to the eigenbasis, i.e.,
\begin{equation}
\bar z = C z,
\label{change variable}
\end{equation}
with column vector $z= (q_1,q_2,p_1,p_2)^T$. We note that $C$ is a symplectic matrix, which requires that the columns be carefully scaled generalized eigenvectors, 
\begin{equation}
\begin{aligned}
C= \begin{pmatrix}
\frac{1}{s_1} & \frac{1}{s_2} & - \frac{1}{s_1} & 0\\
\frac{c_x - \lambda^2}{s_1} & \frac{ \omega_p^2 + c_x}{s_2} & \frac{\lambda^2 - c_x}{s_1}  & 0\\
\frac{\lambda}{s_1} & 0 & \frac{ \lambda}{s_1} & \frac{\omega_p}{s_2} \\
\frac{c_x  \lambda - \lambda^3}{s_1} & 0 & \frac{c_x \lambda - \lambda^3}{s_1} &  \frac{c_x \omega_p + \omega_p^3 }{s_2} 
\end{pmatrix}
\label{sym matrix}
\end{aligned},
\end{equation}
where $s_1$ and $s_2$ are the rescaling factors defined by $s_1=\sqrt{d_{\lambda}}$ and $s_2=\sqrt{d_{\omega_p}}$, respectively, and, 
\begin{equation}
\begin{aligned}
d_\lambda &= \lambda [4 - 2 (c_x -c_y) ( \lambda^2 -c_x)],\\
d_{\omega_p}&= \frac{\omega_p}{2} [4 +2 (c_x - c_y) ( \omega_p^2 + c_x)].
\end{aligned}
\end{equation}
Notice that the symplectic matrix $C$ should satisfy the following relation,
\begin{equation}
C^TJC=J,
\end{equation}
in which $J$ is the $4 \times 4$ canonical symplectic matrix,
\begin{equation}
J=
\begin{pmatrix}
0& I_2\\ -I_2 & 0
\end{pmatrix},
\end{equation}
where $I_2$ is the $2 \times 2$ identity matrix.

Using the new coordinates $q_1$, $q_2$, $p_1$, and $p_2$, the differential equations of the conservative system can be converted to,
\begin{equation}
\begin{aligned}
&\dot q_1=\lambda q_1, \hspace{0.2in} &&\dot p_1=-\lambda p_1,\\
&\dot q_2 = \omega_p, \hspace{0.2in} && \dot p_2=-\omega_p q_2.
\end{aligned}
\label{conservative normal form}
\end{equation}
and the quadratic Hamiltonian function becomes,
\begin{equation}
\mathcal{H}_2=\lambda p_1 q_2 + \tfrac{1}{2} \omega_p \left(q_2^2 + p_2^2\right)
\label{symplectic_Hamiltonian}
\end{equation}
The solutions of \eqref{conservative normal form} can be conveniently written as,
\begin{equation}
\begin{aligned}
& q_1= q_1^0 e^{ \lambda t}, \ \ \ p_1= p_1^0 e^{ - \lambda t},\\
& q_2 + i p_2 = \left(q_2^0 + i  p_2^0 \right) e^{-i \omega_p t}.
\end{aligned}
\label{cons_solut_eigen}
\end{equation}
where the constants $q_1^0$, $p_1^0$, and $q_2^0 + i p_2^0$ are the initial conditions. 
Note that the two functions
\[
f_1 = q_1 p_1, \quad
f_2 =  q_2^2 + p_2^2, 
\]
are constants of motion under the Hamiltonian system \eqref{conservative normal form}, as is $ \mathcal{H}_2$, being a linear combination of $f_1$ and $f_2$. 

\paragraph{The linearized phase space.} 
For positive $h$ and $c$,
the equilibrium or bottleneck region $\mathcal{R}$ (sometimes just called the neck region), which is determined by,
\[
\mathcal{H}_2=h, \quad  \mbox{and} \quad |p_1-q_1|\leq c,
\]
where $c>0$, is homeomorphic to the product of a 2-sphere and an interval $I \in \mathbb{R}$,
$S^2\times I$;
namely, for each fixed value of $p_1 -q_1 $ in the interval $I=[-c,c]$,
we see that  the equation $\mathcal{H}_2=h$ determines a 2-sphere,
\begin{equation}\label{2-sphere}
\tfrac{\lambda }{4}(q_1 +p_1 )^2
+ \tfrac{1}{2}\omega_p (q_2^2+p_2^2)
=h+\tfrac{\lambda }{4}(p_1 -q_1 )^2.
\end{equation}
Suppose $a \in I$, then \eqref{2-sphere} can be re-written as,
\begin{equation}\label{2-sphere2}
x_1^2 + q_2^2+p_2^2
= r^2,
\end{equation}
where $x_1 = \sqrt{\tfrac{1 }{2}\tfrac{\lambda}{\omega_p}}(q_1 +p_1 )$ and 
$r^2=\tfrac{2}{\omega_p}(h+\tfrac{\lambda }{4}a^2)$, which defines a 2-sphere of radius $r$ in the three variables $x_1$, $q_2$, and $p_2$.

\begin{figure}[!t]
	\begin{center}
		\includegraphics[width=\textwidth]{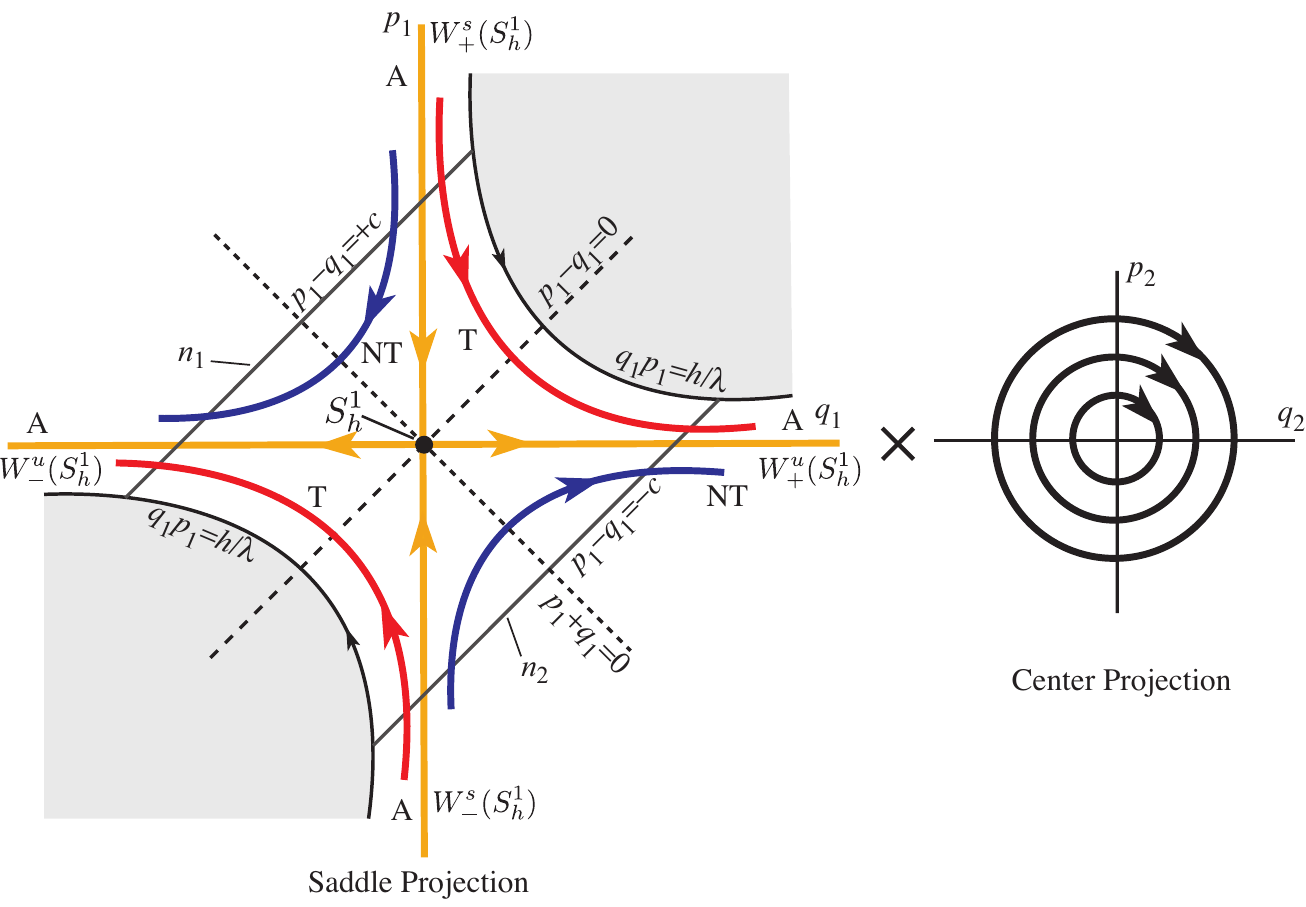}
	\end{center}
	\caption{\footnotesize 
		The flow in the equilibrium region for the conservative system has the form
		saddle $\times$ center.
		On the left is shown a schematic of the projection onto the $(q_1,p_1)$-plane, the saddle projection.
		For the conservative dynamics, the Hamiltonian function $\mathcal{H}_2$ remains constant at $h>0$. Shown are the periodic orbit
		(black dot at the center), the asymptotic orbits (labeled A), two
		transit orbits (T) and two non-transit orbits (NT).
	}
	\label{conservative eigen flow}
\end{figure}

The bounding 2-sphere of $\mathcal{R}$ for which $p_1 -q_1 = c$ will
be called $n_1$ (the ``left'' bounding 2-sphere), and where $p_1 -q_1 = -c$,
$n_2$ (the ``right'' bounding 2-sphere). 
Therefore, $\partial \mathcal{R} =\{ n_1, n_2 \}$.
See Figure \ref{conservative eigen flow}.  

We call the set of points on each
bounding 2-sphere where $q_1 + p_1 = 0$ the equator, and the sets where
$q_1 + p_1 > 0$ or $q_1 + p_1 < 0$ will be called the northern and 
southern hemispheres, respectively.

\paragraph{The linear flow in $\mathcal{R}$.} 
To analyze the flow in
$\mathcal{R}$,   consider
the projections on the ($q_1, p_1$)-plane and the $(q_2,p_2)$-plane, respectively.
In the first case we see the standard
picture of a saddle point in two dimensions,
and in the second, of a center consisting of 
harmonic oscillator motion.
Figure \ref{conservative eigen flow} schematically illustrates the flow.
With regard to the first projection we
see that $\mathcal{R}$ itself projects
to a set bounded on two sides by
the hyperbolas
$q_1p_1 = h/\lambda $
(corresponding to $q_2^2+p_2^2=0$, see \eqref{symplectic_Hamiltonian}) and on two
other sides by the line segments
$p_1-q_1= \pm c$, which correspond to the bounding 2-spheres, $n_1$ and $n_2$, respectively.

Since $q_1p_1$ is an integral
of the equations in $\mathcal{R}$,
the projections of
orbits in the $(q_1,p_1)$-plane
move on the branches of the corresponding
hyperbolas $q_1p_1 =$ constant,
except in the case $q_1p_1=0$, where $q_1 =0$ or $p_1 =0$.
If $q_1p_1 >0$, the branches connect
the bounding line segments $p_1 -q_1 =\pm c$ and if $q_1p_1 <0$, they
have both end points on the same segment.  A check of equation
\eqref{cons_solut_eigen} shows that the orbits move
as indicated by the arrows in Figure \ref{conservative eigen flow}.

To interpret Figure \ref{conservative eigen flow} as a flow in $\mathcal{R}$,  notice
that each point in the $(q_1,p_1)$-plane projection 
corresponds to a 1-sphere, $S^1$, or circle, in
$\mathcal{R}$ given by, 
\[
q_2^2+p_2^2
=\tfrac{2 }{\omega_p}(h-\lambda q_1p_1) .
\]
Of course, for points on the bounding  hyperbolic
segments ($q_1p_1 =h/\lambda $), the
1-sphere collapses to a point. Thus, the segments of the lines 
$p_1-q_1 =\pm c$  in the projection correspond to the 2-spheres
bounding $\mathcal{R}$.  This is because each corresponds to a
1-sphere crossed with an interval where the two end 1-spheres are
pinched to a point.

We distinguish nine classes of orbits grouped into the following four
categories:
\begin{enumerate}
	\item The point $q_1 =p_1 =0$ corresponds to an invariant
	1-sphere $S^1_h$, an unstable {\bf periodic orbit} in
	$\mathcal{R}$ of energy $\mathcal{H}_2=h$.  This 1-sphere is given by,
	\begin{equation}\label{3-sphere}
	q_2^2+p_2^2=\tfrac{2 }{\omega_p}h, 
	\hspace{0.3in} q_1 =p_1 =0.
	\end{equation}
	It is an example of a 
	normally hyperbolic invariant manifold (NHIM) (see \cite{Wiggins1994}).
	Roughly, this means that the stretching and contraction rates under
	the linearized dynamics transverse to the 1-sphere dominate those tangent
	to the 1-sphere.  This is clear for this example since the dynamics normal
	to the 1-sphere are described by the exponential contraction and expansion
	of the saddle point dynamics.  Here the 1-sphere acts as a ``big
	saddle point''.  
	See the black dot at the center of the $(q_1,p_1)$-plane on the left side 
	of Figure
	\ref{conservative eigen flow}.
	
	\item The four half open segments on the axes, $q_1p_1 =0$, 
	correspond to four 
	cylinder surfaces of orbits asymptotic to this invariant 1-sphere 
	$S^1_h$ either as time
	increases ($q_1 =0$) or as time decreases ($p_1 =0$).  These are called {\bf
		asymptotic} orbits and they are the stable and the unstable manifolds of
	$S^1_h$.  The stable manifolds, $W^s_{\pm}(S^1_h)$, are given by,
	\begin{equation}\label{stable_manifold}
	q_2^2+p_2^2=\tfrac{2 }{\omega_p}h, 
	\hspace{0.3in} q_1 =0,
	\hspace{0.3in} p_1 ~{\rm arbitrary}.
	\end{equation}
	$W^s_+(S^1_h)$ (with $p_1>0$) is the branch entering from $n_1$ and
	$W^s_-(S^1_h)$ (with $p_1<0$) is the branch entering from $n_2$.
	The unstable manifolds, $W^u_{\pm}(S^1_h)$, 
	are given by,
	\begin{equation}\label{unstable_manifold}
	q_2^2+p_2^2=\tfrac{2 }{\omega_p}h, 
	\hspace{0.3in} p_1 =0,
	\hspace{0.3in} q_1 ~{\rm arbitrary}
	\end{equation}
	$W^u_+(S^1_h)$ (with $q_1>0$) is the branch exiting from $n_2$ and
	$W^u_-(S^1_h)$ (with $q_1<0$) is the branch exiting from $n_1$.
	See the four orbits labeled A of Figure \ref{conservative eigen flow}.
	
	\item The hyperbolic segments determined by
	$q_1p_1 ={\rm constant}>0$ correspond
	to two solid cylinders of orbits 
	which cross $\mathcal{R}$ from one bounding 2-sphere to the
	other, meeting both in the same hemisphere; the northern hemisphere
	if they go from
	$p_1-q_1 =+c$ to $p_1-q_1 =-c$, and the southern hemisphere
	in the other case. Since
	these orbits transit from one realm to another, we call  them {\bf transit}
	orbits.  See the two orbits labeled T of Figure \ref{conservative eigen flow}.
	
	\item Finally the hyperbolic segments determined by $q_1p_1 = {\rm
		constant}<0$ correspond to two cylinders of orbits in
	$\mathcal{R}$ each of which runs from one hemisphere to the other hemisphere
	on the same bounding 2-sphere.  Thus if $q_1 >0$, the 2-sphere is $n_2$ ($p_1
	-q_1 =-c$) and orbits run from the southern hemisphere
	($q_1 +p_1 <0$) to the northern hemisphere ($q_1
	+p_1 >0$) while the converse holds if $q_1 <0$, where the 
	2-sphere is
	$n_1$. Since these orbits return to the same realm, we call them {\bf
		non-transit} orbits.  See the two orbits labeled NT of Figure \ref{conservative eigen flow}.
\end{enumerate}

We define the transition region, $\mathcal{T}_h$, as the region of initial conditions of a given initial energy $h$ which transit from one side of the neck region to the other.  
This is the set of all transit orbits, which has the geometry of a solid cylinder.   
The transition region, $\mathcal{T}_h$, is made up of one half which goes to the right 
(from $n_1$ to $n_2$),
$\mathcal{T}_{h+}$, defined by $q_1p_1 ={\rm constant}>0$ with both $q_1>0$ and $p_1>0$, and the other half which goes to the left
(from $n_2$ to $n_1$), 
$\mathcal{T}_{h-}$, 
defined by $q_1p_1 ={\rm constant}>0$ with both $q_1<0$ and $p_1<0$.
The boundaries are $\partial \mathcal{T}_{h+}$ and $\partial \mathcal{T}_{h-}$, respectively.
The closure of $\partial \mathcal{T}_{h}$, $\overline{\partial \mathcal{T}_{h}}$, is equal to the boundaries 
$\partial \mathcal{T}_{h+}$ and $\partial \mathcal{T}_{h-}$, along with the periodic orbit $S^1_h$, i.e., 
$\partial \mathcal{T}_{h-} \cup \partial \mathcal{T}_{h+} \cup S^1_h$.

In summary, for the conservative case, the boundary of the transition region, 
$\partial \mathcal{T}_h$, has the topology of a cylinder. The topology of $\partial \mathcal{T}_h$ will be different for the dissipative case, as will be shown in later sections.
For convenience, we may refer to $\partial \mathcal{T}_{h}$ and $\overline{\partial \mathcal{T}_{h}}$ interchangeably.

\begin{figure}
	\begin{center}
		\includegraphics[width=\textwidth]{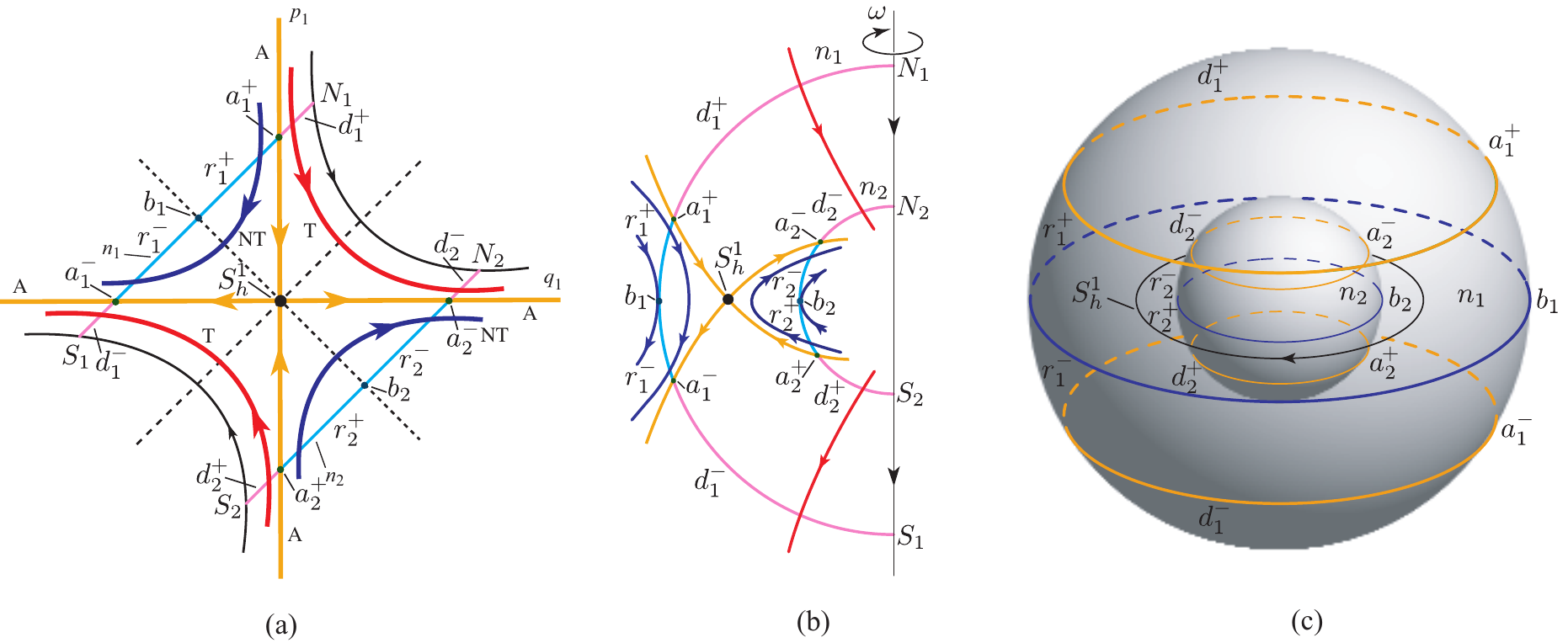}
	\end{center}
	\caption{\footnotesize 
		(a) The projection onto the $(q_1,p_1)$-plane, the saddle projection, with labels consistent with the text and (b) and (c).
		(b) The cross-section of the flow in the
		$\mathcal{R}$ region of the energy surface. The north and south poles of bounding sphere $n_i$ are labeled as $N_i$ and $S_i$, respectively.
		(c) The McGehee representation of the flow on the boundaries of the $\mathcal{R}$ region, highlighting the features on the bounding spheres $n_1$ and $n_2$ for $h>0$.
	}
	\label{conservative McGehee}
\end{figure}

\paragraph{McGehee representation.}
McGehee \cite{McGehee1969}, building on the work of Conley  \cite{Conley1968}, proposed a
representation which makes it easier to visualize the region $\mathcal{R}$.
Recall that $\mathcal{R}$ is a 3-dimensional manifold that is homeomorphic to $S^2\times I$. In \cite{McGehee1969},
it is represented by a spherical annulus  bounded by the 
two 2-spheres $n_1, n_2$, as shown in Figure \ref{conservative McGehee}(c).

Figure \ref{conservative McGehee}(a) is a cross-section of $\mathcal{R}$. Notice that this
cross-section is qualitatively the same as the saddle projection illustration in Figure
\ref{conservative eigen flow}.  
The full picture (Figure \ref{conservative McGehee}(c)) is obtained by rotating
this cross section, Figure \ref{conservative McGehee}(b), about the indicated axis, where the azimuthal angle $\omega$ roughly describes the angle in the center projection in Figure \ref{conservative eigen flow}.
The following
classifications of orbits correspond to the previous four categories:
\begin{enumerate}
	\item  There is an invariant
	1-sphere $S^1_h$, a {\it periodic orbit} 
	in the region $\mathcal{R}$ corresponding to the black dot in the
	middle of Figure \ref{conservative McGehee}(a).
	Notice that this 1-sphere is the equator of the
	central 2-sphere given by $p_1 -q_1 =0$. 
	
	\item Again let $n_1,n_2$ be the bounding 2-spheres
	of region $\mathcal{R}$, and let $n$ denote either $n_1$ or $n_2$.
	We can divide $n$ into two hemispheres: $n^+$, where the flow enters
	$\mathcal{R}$,
	and $n^-$, where the flow leaves $\mathcal{R}$.  
	There are four cylinders of orbits 
	asymptotic to the invariant 1-sphere $S^1_h$.
	They form the stable and unstable manifolds which are {\it asymptotic} to the invariant 1-sphere
	$S^1_h$.  
	Topologically, both invariant manifolds 
	look like 2-dimensional cylinders or ``tubes'' ($S^1\times {\mathbb R}$) inside
	a 3-dimensional energy manifold.
	The interior of the stable manifolds  
	$W^s_{\pm}(S^1_h)$ and unstable manifolds
	$W^u_{\pm}(S^1_h)$ can be given as follows 
	\begin{equation}\label{interior}
	\begin{split}
	{\rm int}(W^s_+(S^1_h))
	&=
	\{(q_1,p_1,q_2,p_2)\in \mathcal{R}\mid 
	\hspace{.1in} p_1>q_1>0\},  \\
	{\rm int}(W^s_-(S^1_h))
	&=
	\{(q_1,p_1,q_2,p_2)\in \mathcal{R} \mid 
	\hspace{.1in} p_1<q_1<0\},  \\
	{\rm int}(W^u_+(S^1_h))
	&=
	\{(q_1,p_1,q_2,p_2)\in \mathcal{R}\mid 
	\hspace{.1in} q_1>p_1>0\},  \\
	{\rm int}(W^u_-(S^1_h))
	&=
	\{(q_1,p_1,q_2,p_2)\in \mathcal{R}\mid 
	\hspace{.1in} q_1<p_1<0\}.
	\end{split}
	\end{equation}
	The exterior of these invariant manifolds can be given similarly from
	studying Figure \ref{conservative McGehee}(a) and (b).
	
	\item
	Let $a^+$ and $a^-$
	(where $q_1 =0$ and $p_1 =0$ respectively) be the
	intersections of the stable and unstable manifolds with 
	the bounding sphere $n$.
	Then $a^+$ appears as a 1-sphere in $n^+$, and $a^-$ appears as a 
	1-sphere in $n^-$.
	Consider the two spherical caps on each bounding 
	2-sphere given by
	\begin{align*}
	d_1^+&=\{(q_1,p_1,q_2,p_2)\in \mathcal{R}
	\mid \hspace{.1in}p_1-q_1 =+c, 
	\hspace{.1in}p_1>q_1 >0\},\\
	d_1^-&=\{(q_1,p_1,q_2,p_2)\in \mathcal{R}
	\mid \hspace{.1in}p_1 -q_1 =+c,
	\hspace{.1in}q_1<p_1<0\},\\
	d_2^+ &= \{(q_1,p_1,q_2,p_2)\in \mathcal{R}\mid 
	\hspace{.1in}p_1 -q_1 =-c, \hspace{.1in}
	p_1<q_1 <0\}, \\
	d_2^-&=\{(q_1,p_1,q_2,p_2)\in \mathcal{R}
	\mid \hspace{.1in}p_1 -q_1 =-c,
	\hspace{.1in} q_1 >p_1>0\}.
	\end{align*}
	Since $d_1^+$ is  the spherical cap 
	in $n_1^+$ bounded by $a_1^+$, then
	the {\it transit} orbits
	entering $\mathcal{R}$ on $d_1^+$ exit on $d_2^-$ of the other bounding sphere.
	Similarly, since 
	$d_1^-$ is  the spherical cap in $n_1^-$ bounded by $a_1^-$, the transit
	orbits leaving on $d_1^-$ have
	come from $d_2^+$ on the other bounding sphere.
	Note that all spherical caps where the transit orbits pass through are
	in the interior of stable and unstable manifold tubes. 
	
	\item Let $b$ be the intersection
	of $n^+$ and $n^-$ (where $q_1 +p_1 =0$).  Then, $b$
	is a 1-sphere of tangency
	points.  Orbits tangent at this 1-sphere ``bounce off,'' i.e., do not
	enter
	$\mathcal{R}$ locally.  Moreover, if we let $r^+$ be a spherical zone which
	is bounded by $a^+$ and $b$, then {\it non-transit} orbits
	entering $\mathcal{R}$ on $r^+$ 
	exit on the same bounding 2-sphere through $r^-$ 
	which is bounded by
	$a^-$ and $b$.  It is easy to show that all the spherical zones where  
	non-transit orbits bounce off are in the exterior of stable and unstable
	manifold tubes.
\end{enumerate}

The McGehee representation provides an additional, perhaps clearer, visualization of the dynamics in the equilibrium region. In particular, the features on the two spheres, $n_1$ and $n_2$, which form $\partial \mathcal{R}$ for a constant $h>0$, can be considered in the dissipative case as well, and compared with the situation in the conservative case, as shown for some examples below.  The spheres $n_1$ and $n_2$ can be viewed as spherical Poincar\'e sections parametrized by their distance from the saddle point, $c$, which reveal the topology of the transition region boundary, $\partial \mathcal{T}_h$, particularly through how the geometry of $a_i^+$ and $a_i^-$ (for $i=1,2$) change as $c$ changes.


\subsection{Dissipative system}
\label{linear_dissipative_system}
For the dissipative system, we still use the symplectic matrix $C$ in \eqref{sym matrix} to perform a transformation, via \eqref{change variable}, to the symplectic eigenspace, even though this is no longer the true eigenspace of the dissipative linearization matrix $A=M+D$. The equations of motion in the symplectic eigenspace  are, 
\begin{equation}
\dot z = \Lambda z + \Delta z,
\label{B:EOM with damping in phase space}
\end{equation}
where $\Lambda=C^{-1}MC$ is the conservative part of the dynamics, as before, and the transformed damping matrix is,
\begin{equation}
\Delta= C^{-1}DC= -c_h
\begin{pmatrix}
\tfrac{1}{2} & 0 & \tfrac{1}{2} & 0\\
0 & 0 & 0 & 0\\
\tfrac{1}{2} & 0 & \tfrac{1}{2} & 0\\
0 & 0 & 0 & 1
\end{pmatrix}.\label{standard_damping_matrix}
\end{equation}

To analyze the behavior in the dissipative eigenspace (as opposed to the symplectic eigenspace), 
the eigenvalues and eigenvectors,
$\beta_i$ and $u_{\beta_i}$, respectively, $(i=1,...,4)$, are,
\begin{equation}
\begin{aligned}
\beta_{1,2} &=-\delta \mp \tfrac{1}{2}\sqrt{c_h^2 + 4 \lambda^2}, 
\hspace{0.2in} && u_{\beta_{1,2}}=\left(\delta,0, \lambda \pm \tfrac{1}{2} \sqrt{c_h^2 + 4 \lambda^2},0 \right)^T,\\
\beta_{3,4}& =-\delta \pm i \omega_d, \hspace{0.2in} &&u_{\beta_{3,4}}=\left(0, \omega_p, 0, - \delta \pm i \omega_d \right)^T,
\end{aligned}
\end{equation}
where $\delta=\tfrac{1}{2}c_h$, $\omega_d=\omega_p \sqrt{1-\xi_d^2}$ and $\xi_d=\delta/\omega_p$.
Thus, the general (real) solutions are,
\begin{equation}
\begin{aligned}
& q_1(t) =k_1 e^{\beta_1 t} + k_2 e^{\beta_2 t}, 
\hspace{0.2in} 
p_1(t)=k_3 e^{\beta_1 t} + k_4 e^{\beta_2 t},\\
& q_2(t) = k_5  e^{- \delta t} \cos{\omega_d t} + k_6 e^{- \delta t} \sin{\omega_d t},\\
& p_2(t) = \frac{k_5 }{\omega_p} e^{- \delta t} \left(-\delta \cos{\omega_d t} - \omega_d \sin{\omega_d t} \right) +\frac{k_6 }{\omega_p}  e^{- \delta t} \left(\omega_d \cos{\omega_d t - \delta \sin{\omega_d t}} \right),\\
\end{aligned}
\label{dissip_solut_eigen}
\end{equation}
where,
\begin{equation*}
\begin{aligned}
k_1 &= \frac{q_1^0 \left(2 \lambda + \sqrt{c_1^2 + 4 \lambda^2} \right)-c_1 p_1^0 }{2\sqrt{c_1^2 + 4 \lambda^2}}, \hspace{0.2in}&& k_2 =\frac{q_1^0 \left(-2 \lambda + \sqrt{c_1^2 + 4 \lambda^2} \right)+c_1 p_1^0 }{2\sqrt{c_1^2 + 4 \lambda^2}},\\
k_3 &= \frac{p_1^0 \left(-2 \lambda + \sqrt{c_1^2 + 4 \lambda^2} \right)-c_1 q_1^0 }{2\sqrt{c_1^2 + 4 \lambda^2}}, && k_4 = \frac{p_1^0 \left(2 \lambda + \sqrt{c_1^2 + 4 \lambda^2} \right)+c_1 q_1^0 }{2\sqrt{c_1^2 + 4 \lambda^2}},\\
k_5&=q^0_2 , \hspace{0.5in} k_6=\frac{p^0_2 \omega_p + q^0_2 \delta}{\omega_d}.
\end{aligned}
\end{equation*}
Taking the total derivative of the Hamiltonian with respective to time along trajectories and using \eqref{B:EOM with damping in phase space}, we have,
\begin{equation*}
\frac{\mathrm{d} \mathcal{H}_2}{\mathrm{d} t}= - \tfrac{1}{2} c_h \lambda \left(q_1 + p_1 \right)^2 - c_h \omega_p p_2^2 \le 0,
\end{equation*}
which means the Hamiltonian is generally decreasing (more precisely, non-increasing) due to damping. 

\paragraph{The linear flow in $\mathcal{R}$.} 
Similar to the discussions in the conservative system, we still choose the same equilibrium region $\mathcal{R}$ to consider the projections on the $\left(q_1, p_1\right)$-plane and $\left(q_2, p_2\right)$-plane, respectively. Different from the saddle $\times$ center projections in the conservative system, here we see saddle $\times$ focus projections in the dissipative system. The stable focus is a damped oscillator with frequency of $\omega_d=\omega_p \sqrt{1-\xi_d^2}$. Different classes of orbits can also be grouped into the following four categories:
\begin{enumerate}
	\item The point $q_1=p_1=0$ corresponds to a  {\bf focus-type asymptotic} orbit with motion purely in the $(q_2,p_2)$-plane (see black dot at the origin of
	the $(q_1,p_1)$-plane in Figure \ref{B: dissipative eigen flow}).  
	\begin{figure}[!t]
		\begin{center}
			\includegraphics[width=\textwidth]{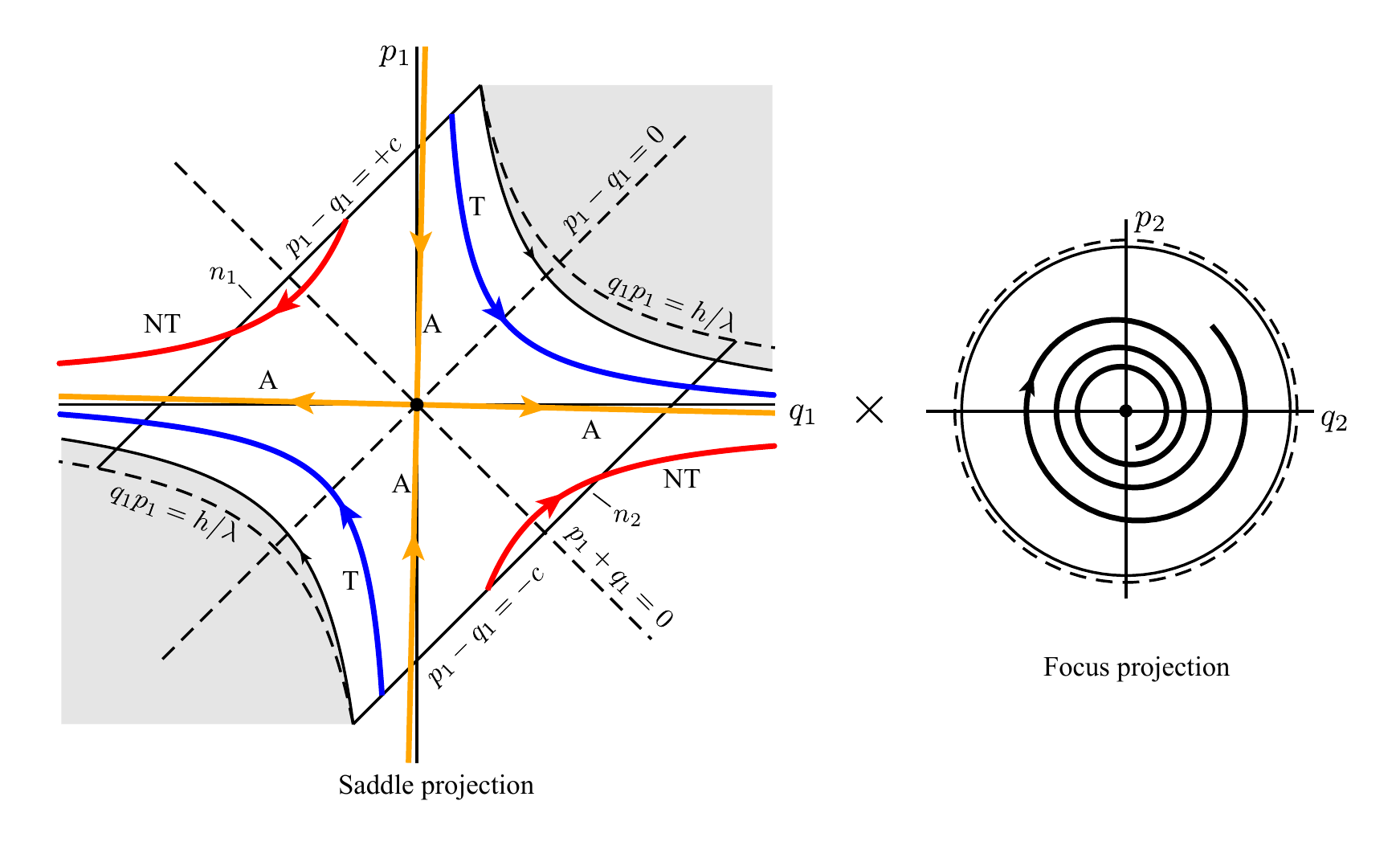}
		\end{center}
		\caption{\footnotesize 
			The flow in the equilibrium region for the dissipative system has the form saddle $\times$ focus. On the left is shown the saddle projection onto the $(q_1,p_1)$-plane. The black dot at the origin represents focus-type asymptotic orbits with only a focus projection, thus oscillatory dynamics decaying towards the equilibrium point. The asymptotic orbits (labeled A) are the saddle-type asymptotic orbits which are tilted clockwise compared to the conservative system. They still form the separatrix between transit orbits (T) and non-transit orbits (NT). The hyperbolas, $q_1 p_1=h/\lambda$, are no longer the boundary of trajectories with initial conditions on the bounding sphere ($n_1$ or $n_2$) due to the dissipation of the energy. The boundary of the shaded region are still the fastest trajectories with initial conditions on the bounding sphere, but are not strictly hyperbolas. Note that the saddle projection and focus projection are uncoupled in this dissipative system.
		}
		\label{B: dissipative eigen flow}
	\end{figure}
	Such orbits are asymptotic to the equilibrium point itself, rather than a periodic orbit of energy $h$ as in the conservative case.
	Due to the effect of damping, the periodic orbits on each energy manifold of energy $h$ do not exist.
	The 1-sphere $S_h^1$ still exists, but is no longer invariant. Instead, it corresponds to all the initial conditions of initial energy $h$ which are focus-type asymptotic orbits.  
	The projection of $S_h^1$ to the configuration space in the dissipative system is the same as the projection of the periodic orbit in the conservative system.
	
	\item The four half open segments on the lines governed by $q_1=  c_h p_1/(2 \lambda \pm \sqrt{c_1^2 + 4 \lambda^2}) $ correspond to {\bf saddle-type asymptotic} orbits.
	See the four orbits labeled A in Figure \ref{B: dissipative eigen flow}. 
	
	\item The segments which cross $\mathcal{R}$ from one boundary to the other, i.e., from $p_1 - q_1=+c$ to $p_1 - q_1=-c$ in the northern hemisphere, and vice versa in the southern hemisphere, correspond to {\it transit} orbits. See the two orbits labeled $T$ of Figure \ref{B: dissipative eigen flow}.
	
	\item  Finally the segments which run from one hemisphere to the other hemisphere on the same boundary, namely which start from $p_1 - q_1 = \pm c$ and return to the same boundary, correspond to {\it non-transit} orbits. See the two orbits labeled NT of Figure \ref{B: dissipative eigen flow}.
\end{enumerate}

As done in Section \ref{linear_conservative_system}, we define the transition region, $\mathcal{T}_h$, as the region of initial conditions of a given initial energy $h$ which transit from one side of the neck region to the other.  
As before, the transition region, $\mathcal{T}_h$, is made up of one half which goes to the right,
$\mathcal{T}_{h+}$, 
and the other half which goes to the left, 
$\mathcal{T}_{h-}$. 
The boundaries are $\partial \mathcal{T}_{h+}$ and $\partial \mathcal{T}_{h-}$, respectively.
The closure of $\partial \mathcal{T}_{h}$, $\overline{\partial \mathcal{T}_{h}}$, is equal to the boundaries 
$\partial \mathcal{T}_{h+}$ and $\partial \mathcal{T}_{h-}$, along with the focus-type asymptotic initial conditions $S^1_h$, i.e., as before,
$\partial \mathcal{T}_{h-} \cup \partial \mathcal{T}_{h+} \cup S^1_h$.

As shown below, for the dissipative case, the closure of the boundary of the transition region, 
$\partial \mathcal{T}_h$, has the topology of an ellipsoid, rather than a cylinder as in the conservative case.  
As before, for convenience, we may refer to $\partial \mathcal{T}_{h}$ and $\overline{\partial \mathcal{T}_{h}}$ interchangeably.

\paragraph{McGehee representation.}
Similar to the McGehee representation for the conservative system given in Section \ref{linear_conservative_system} to visualize the region $\mathcal{R}$, here we utilize the McGehee representation again to illustrate the behavior in same region for the dissipative system. 
All labels are consistent throughout the paper.

\begin{figure}
	\begin{center}
		\includegraphics[width=\textwidth]{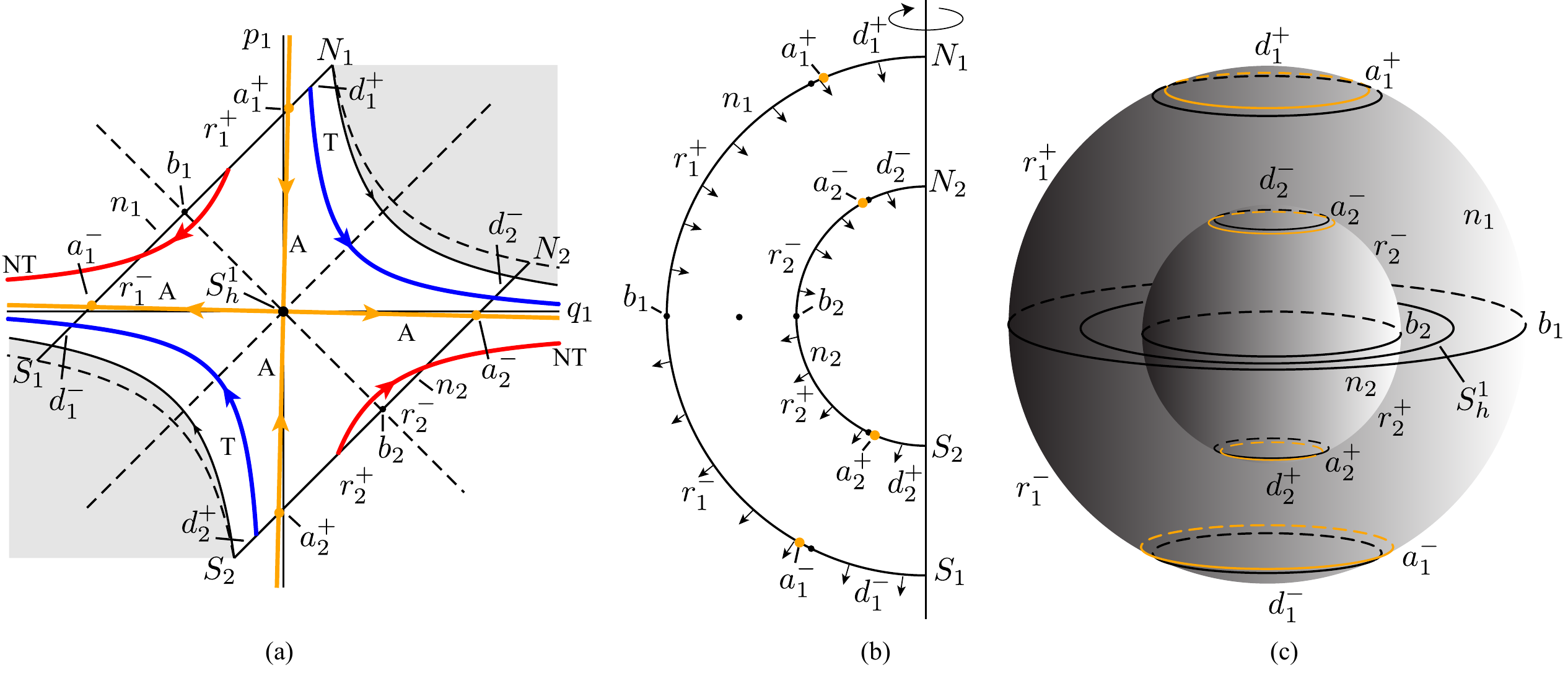}
	\end{center}
	\caption{\footnotesize 
		(a) The projection onto the $(q_1,p_1)$-plane, the saddle projection, with labels consistent with the text and (b) and (c).
		(b) The cross-section of the flow in the
		$\mathcal{R}$ region of the energy surface. The north and south poles of bounding sphere $n_i$ are labeled as $N_i$ and $S_i$, respectively.
		(c) The McGehee representation of the
		flow in the region $\mathcal{R}$.
	}
	\label{B: McGehee representation}
\end{figure}

Note that since the McGehee representation uses spheres with the same energy to show the dynamical behavior in phase space, while the energy of any particular trajectory in the dissipative system decreases gradually during evolution, Figures \ref{B: McGehee representation}(b) and \ref{B: McGehee representation}(c) show only the initial conditions at a given initial energy. Therefore, in the present McGehee representation, only the initial conditions on the two bounding spheres are shown and discussed in the next part. In addition, the black dot near the orange dots $a_i^{\pm}$ and $b_i^{\pm}$ ($i=1,2$) in Figure \ref{B: McGehee representation}(b) are the corresponding dots in the conservative system which are used to show how damping affects the transition. 

The following classifications of orbits correspond to the previous four categories:
\begin{enumerate}
	\item 1-sphere $S_h^1$ exists in the region $\mathcal{R}$ corresponding to the black dot in the middle of Figure \ref{B: McGehee representation}(b) and the equator of the central 2-sphere given by $p_1-q_1=0$ in \ref{B: McGehee representation}(c). The 1-sphere gives the initial conditions of the initial energy $h$ for all focus-type asymptotic orbits. The same 1-sphere in the conservative system is invariant under the flow, that is, a periodic orbit of constant energy $h$. However, the corresponding $S_h^1$ is not invariant in the dissipative system, since the energy is decreasing during evolution due to the damping.

	\item There are four 1-spheres in the region $\mathcal{R}$ starting in the bounding 2-spheres $n_1$ and $n_2$ which give the initial conditions for orbits asymptotic to the equilibrium point. Two of them in $n^+$, labeled by $a^+$, are stable saddle-type asymptotic orbits and the other two in $n^-$, labeled by $a^-$, are unstable  asymptotic orbits, where $a^+$ and $a^-$ are given by,
	\begin{equation}
	\begin{aligned}
	&a_1^+=\left\{(q_1,p_1,q_2,p_2)\in \mathcal{R} | \hspace{0.1in} (q_1,p_1)=(k_p,1)c/(1-k_p)\right\},\\
	&a_1^-=\left\{(q_1,p_1,q_2,p_2)\in \mathcal{R} | \hspace{0.1in} (q_1,p_1)=(-1,k_p)c/(1+k_p)\right\},\\
	&a_2^+=\left\{(q_1,p_1,q_2,p_2)\in \mathcal{R} | \hspace{0.1in} (q_1,p_1)=(k_p,1)c/(k_p-1)\right\},\\
	&a_2^-=\left\{(q_1,p_1,q_2,p_2)\in \mathcal{R} | \hspace{0.1in} (q_1,p_1)=(1,-k_p)c/(1+k_p)\right\},\\
	\end{aligned}
	\end{equation}
	where $k_p = c_h /(2 \lambda + \sqrt{c_h^2 + 4 \lambda^2} )$. As shown in Figure \ref{B: McGehee representation}(c), $a^+$ appears as an orange circle in $n^+$, and $a^-$ appears as an orange circle in $n^-$. The corresponding curves for the same energy in the conservative system are shown as black curves.
	
	\item Consider the two spherical caps on each bounding 2-sphere, $n_1$ and $n_2$, given by,
	\begin{equation}
	\begin{aligned}
	& d_1^+=\left\{(q_1,p_1,q_2,p_2) \in \mathcal{R} \mid  \hspace{0.1in}p_1-q_1=c,  \hspace{0.1in} q_1>ck_p/(1-k_p)\right\},\\
	& d_1^-=\left\{(q_1,p_1,q_2,p_2) \in \mathcal{R} \mid  \hspace{0.1in}p_1-q_1=c,  \hspace{0.1in} q_1<-c/(1+k_q)\right\},\\
	& d_2^+=\left\{(q_1,p_1,q_2,p_2) \in \mathcal{R} \mid  \hspace{0.1in}p_1-q_1=-c,  \hspace{0.1in} q_1<ck_p/(k_p-1)\right\},\\
	& d_2^-=\left\{(q_1,p_1,q_2,p_2) \in \mathcal{R} \mid  \hspace{0.1in}p_1-q_1=-c,  \hspace{0.1in} q_1>c/(1+k_p)\right\}.
	\end{aligned}
	\end{equation}
	The spherical cap $d_1^+$, bounded by the $a_1^+$ on $n_1^+$, gives all initial conditions of initial energy $h$ for the transit orbits starting from the bounding sphere $n_1^+$ and entering $\mathcal{R}$. Similarly, the spherical cap $b_1^-$ in $n_1^-$, bounded by $a_1^-$, determines all initial conditions of initial energy $h$ for transit orbits starting on the bounding sphere $n_1^-$ and leaving $\mathcal{R}$. The spherical caps $d_2^+$ and $d_2^-$ on $n_2$ have similar dynamical behavior. Note that in the conservative system the transit orbits entering $\mathcal{R}$ on $d^+$ will leave on $d^-$ in the same 2-sphere. However, those transit orbits with the same initial conditions in the dissipative system will not leave on the corresponding 2-sphere, but leave on another sphere with lower energy. Moreover, the spherical caps $d^+$ shrink and $d^-$ expand compared to that of the conservative system. Since the area of the caps $d^+$ and $b^-$ determines the amount of transit orbits and non-transit orbits respectively, the shrinkage of the caps $d^+$ and expansion of the caps $d^-$ means the damping reduces the probability of transition and increase the probability of non-transition, respectively. 
	
	\item Let $b$ be the intersection of $n^+$ and $n^-$ (where $q_1+p_1=0$). Then, $b$ is a 1-sphere of tangency points. Orbits tangent at this 1-sphere ``bounce off'', i.e., do not enter $\mathcal{R}$ locally. 
	The spherical zones $r_1$ and $r_2$, bounded by $a^+_i$ and $a^-_i$, give the initial conditions for non-transit orbits zone. $r^+$, bounded by $a^+_i$ and $b_i$, are the initial conditions of initial energy $h$ for non-transit orbits entering $\mathcal{R}$ and $r^-_i$ are the initial conditions of initial energy $h$ for non-transit orbits leaving $\mathcal{R}$. Note that unlike the shift of the spherical caps in the dissipative system compared to that of the conservative system, the tangent spheres $b_1$ and $b_2$ do not move when damping is taken into account. Moreover, in the conservative system, non-transit orbits enter $\mathcal{R}$ on $r^+$ and then exit on the same energy bounding 2-sphere through $r^-$, but the non-transit orbits in the dissipative system exit on a different 2-sphere with different energy determined by the damping and the initial conditions.
\end{enumerate}

\subsection{Transition tube and transition ellipsoid}
\label{B: tube and ellipsoid}

After examining the flow in the eigenspace, we study the linearized dynamics in the phase space.

\paragraph{Transition boundary in the conservative system.}
Once we obtained the analytical solutions in \eqref{cons_solut_eigen} for equations \eqref{conservative normal form} in the eigenspace of the conservative system, we can use the change of variables in \eqref{change variable} to get the analytical solutions for the equations in the conservative system written as
\begin{equation}
\begin{aligned}
& \bar q_1 = \frac{q_1^0}{s_1} e^{\lambda t} - \frac{p_1^0}{s_1} e^{-\lambda t}  + \frac{1}{s_2} (q_2^0 \cos \omega_p t + p_2^0 \sin \omega_p t),\\
& \bar q_2=\frac{c_x -\lambda^2 }{s_1} q_1^0 e^{\lambda t} + \frac{\lambda^2 - c_x}{s_1} p_1^0 e^{-\lambda t} + \frac{\omega_p^2 + c_x }{s_2} (q_2^0 \cos \omega_p t + p_2^0 \sin \omega_p t),\\
& \bar p_1 = \frac{\lambda}{s_1} q_1^0 e^{\lambda t} + \frac{\lambda}{s_1} p_1^0 e^{-\lambda t} + \frac{\omega_p}{s_2} (p_2^0 \cos \omega_p t- q_2^0 \sin \omega_p t),\\
& \bar p_2 = \frac{c_x \lambda - \lambda^3}{s_1} q_1^0 e^{\lambda t} + \frac{c_x \lambda - \lambda^3}{s_1} p_1^0 e^{-\lambda t} + \frac{c_x \omega_p + \omega_p^3}{s_2}(p_2^0 \cos \omega_p t- q_2^0 \sin \omega_p t).
\end{aligned}
\label{conserv_solution_phase}
\end{equation}
From the discussion about the conservative system in Section \ref{linear_conservative_system}, we know that the invariant manifold of the periodic orbit acts as a separatrix which separate two distinct types of motion: transit orbits and non-transit orbits. Thus, we can compute the initial conditions of the asymptotic orbits to get the transition boundary of a given energy $h$. To get such initial conditions, we need to set the coefficient of the unstable term $e^{\lambda t}$ in $\bar q_1$ as zero, since this term will go to infinity along positive time which is against the asymptotic properties. Thus, we have
\begin{equation}
q_1^0=0.
\end{equation}
Denoting the initial conditions in the phase space by $\bar q_1^0$, $\bar q_2^0$, $\bar p_1^0$, and $\bar p_2^0$, they can be connected to the initial conditions in the eigenspace via \eqref{conserv_solution_phase} by taking $t=0$ from which we can straightforwardly write $q_1^0$, $q_2^0$, $p_1^0$, $p_2^0$, and $\bar p_1^0$ in terms of $\bar q_2^0$, $\bar p_1^0$, and $\bar p_2^0$. In this case, the normal form of Hamiltonian function in \eqref{symplectic_Hamiltonian} can be rewritten by
\begin{equation}
\frac{\left[(\lambda^2 - c_x)x +y\right]^2}{\frac{2h(\lambda^2+\omega_p^2)^2}{s_2^2 \omega_p}} + \frac{\left[-\lambda x +(\lambda^2 - c_x)\lambda y +(\lambda^2 + \omega_p^2)p_y \right]^2}{\frac{2h \omega_p (c_x + \omega_p^2)^2 (\lambda^2 + \omega_p^2)^2}{s_2^2}}=1
\label{tube_formula}
\end{equation}

Notice the formula in \eqref{tube_formula} is a tube.  Its projection onto the position space is a strip bounded by the lines \cite{zhong2018tube},
\begin{equation}
\bar q_2 = (c_x - \lambda^2) \bar q_1 \pm  \frac{\lambda^2 + \omega_p^2}{s_2} \sqrt{\frac{2h}{\omega_p}}.
\label{strip}
\end{equation}

\paragraph{Transition boundary in the dissipative system.}
After we get the solution of the dissipative system in the eigenspace in \eqref{dissip_solut_eigen}, we can use the change of coordinates in \eqref{change variable} to obtain the solutions of the equations \eqref{nond eq} as,
\begin{equation}
\begin{aligned}
&\bar q_1= \frac{k_1 - k_3}{s_1} e^{\beta_1 t} + \frac{k_2-k_4}{s_1} e^{\beta_2 t}+ \frac{q_2}{s_2},\\
& \bar q_2= \frac{k_1-k_3}{s_1} (c_x - \lambda^2) e^{\beta_1 t} + \frac{k_2-k_4}{s_1} (c_x - \lambda^2) e^{\beta_2 t} + \frac{\omega_p^2 + c_x}{s_2} q_2,\\
& \bar p_1= \frac{k_1 + k_3}{s_1} \lambda e^{\beta_1 t} +\frac{k_2 + k_4}{s_1} \lambda e^{\beta_2 t} + \frac{\omega_p}{s_2}p_2,\\
& \bar p_2 = \frac{k_1 + k_3}{s_1} \left(c_x \lambda - \lambda^3 \right)e^{\beta_1 t} +\frac{k_2 + k_4}{s_1} \left(c_x \lambda - \lambda^3 \right)e^{\beta_2 t} + \frac{c_x \omega_p + \omega_p^3}{s_2} p_2.
\end{aligned}
\label{dissi_solution_phase}
\end{equation}
To obtain the initial conditions for the asymptotic orbits in the dissipative system, the term of $e^{\beta_1 t}$ should vanish, otherwise the trajectory will go to infinity along positive time. Thus, we can directly set the coefficient of $e^{\beta_1 t}$ as zero, i.e., $k_1-k_3=0$, which results in
\begin{equation}
q_1^0=k_p p_1^0
\label{dissi_stable_asymp}
\end{equation}
which are the initial conditions of the stable asymptotic orbits in the eigenspace of the dissipative system as discussed before Section \ref{linear_dissipative_system}. Denoting the initial conditions in the phase space by $\bar q_1^0$, $\bar q_2^0$, $\bar p_1^0$, and $\bar p_2^0$, they can be connected to the initial conditions in the eigenspace via \eqref{dissi_solution_phase} by taking $t=0$. Thus, along with \eqref{dissi_stable_asymp}, we can straightforwardly write $q_1^0$, $q_2^0$, $p_1^0$, and $p_2^0$ 
in terms of $\bar q_1^0$, $\bar q_2^0$, $\bar p_1^0$, and $\bar p_2^0$. 
In this case,  the normal form of Hamiltonian function in \eqref{symplectic_Hamiltonian} can be rewritten by
\begin{equation}
\frac{\left[(\lambda^2 - c_x)x +y\right]^2}{\frac{2h(\lambda^2+\omega_p^2)^2}{s_2^2 \omega_p}}+\frac{\left[ (c_x + \omega_p^2)x - y\right]^2}{\frac{h (k_p -1)^2(\lambda^2 + \omega_p^2)^2}{\lambda k_p s_1^2}} + \frac{\left[-\lambda x +  (\lambda^2 - c_x) \lambda  y + \frac{(1-k_p)}{(1+k_p)}(\lambda^2 + \omega_p^2) p_y\right]^2}{\frac{2 h \omega_p(k_p-1)^2(c_x + \omega_p^2)^2 (\lambda^2 + \omega_p^2)^2}{s_2^2 (1+k_p)^2}}=1,
\label{ellipsoid_fourmula}
\end{equation}
which has the form of an ellipsoid. It can also be written in the form, $a_{\bar p_2} \left( \bar p_2^0 \right)^2 + b_{\bar p_2} \bar p_2^0 + c_{\bar p_2}=0$, where $a_{\bar p_2}$, $b_{\bar p_2}$, and $c_{\bar p_2}$ are given in Appendix \ref{Appendix_quadratic}. The projection of the ellipsoid onto the configuration space can be obtained by $b_{\bar p_2}^2 - 4 a_{\bar p_2} c_{\bar p_2}=0$ which is an ellipse of the following form \cite{zhong2018tube},
\begin{equation}
\frac{(\bar q_1^0 \cos \theta + \bar q_2^0 \sin \theta)^2}{a_e^2} + \frac{(-\bar q_1^0 \sin \theta + \bar q_2^0 \cos \theta)^2}{b_e^2}=1,
\label{ellipse}
\end{equation}
where
\begin{equation}
\begin{split}
& 
a_e= \sqrt{\frac{2 h \left(\lambda^2 + \omega_p^2\right)^2 \left(c_x + \omega_p^2 \right)^2}{\omega_p s_2^2 \left[ \left(c_x + \omega_p^2 \right)^2 + 1\right]}}, \hspace{0.7in} 
b_e= \sqrt{\frac{h \left(k_p - 1\right)^2 \left(\lambda^2 + \omega_p^2 \right)^2}{\lambda k_p s_1^2 \left[ \left(c_x + \omega_p^2 \right)^2 + 1\right]}},\\
&  
\cos \theta = \frac{1}{\sqrt{\left(c_x + \omega_p^2 \right)^2 + 1 }}, \hspace{1.3in}   
\sin \theta= \frac{\left(c_x + \omega_p^2 \right)}{\sqrt{\left(c_x + \omega_p^2 \right)^2 + 1 } }\\
\end{split}
\end{equation}

\begin{figure}[!t]
	\begin{center}
		\includegraphics[width=\textwidth]{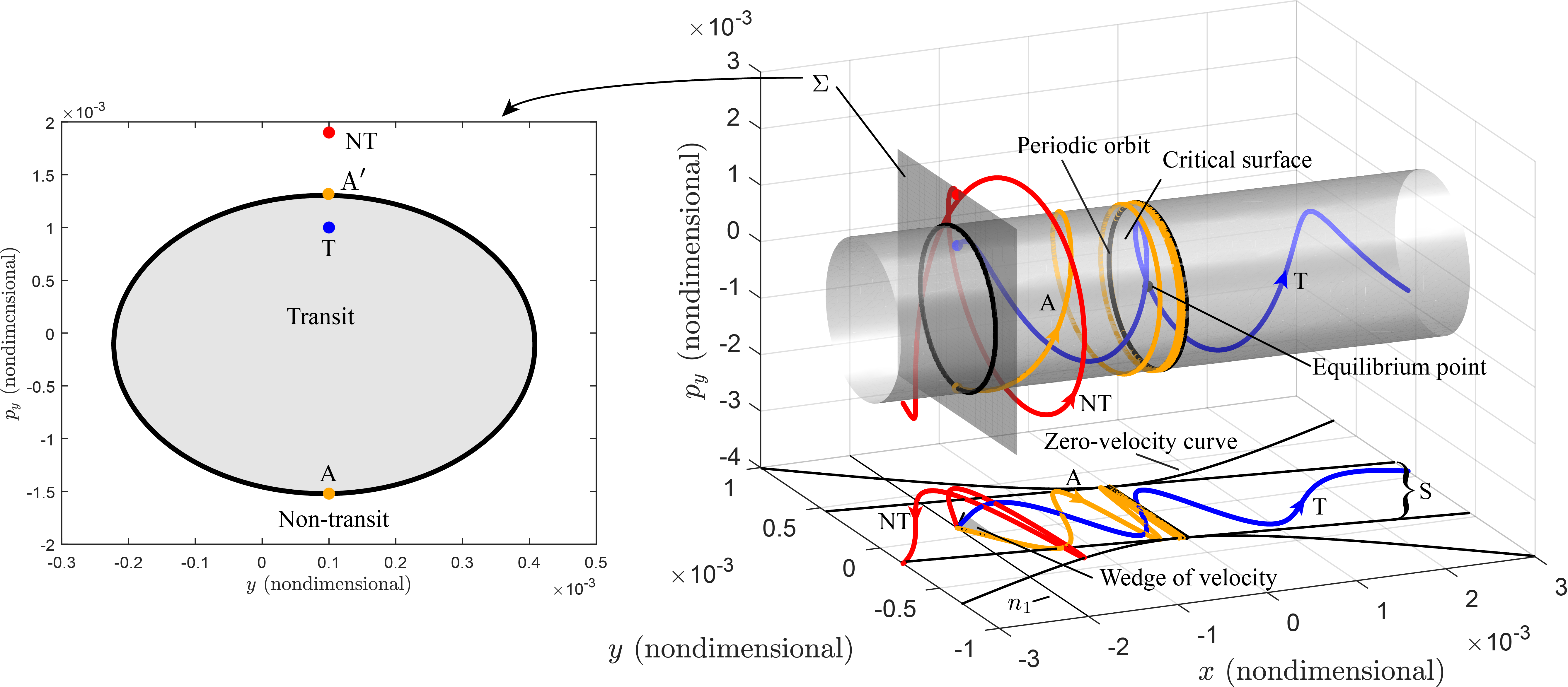}
	\end{center}
	\caption{\footnotesize 
		Transition region boundary $\partial \mathcal{T}_h$ which is a tube (cylinder) for the linearized conservative system with initial energy $h$. The left figure shows tube boundary (the ellipse) separating the transit and non-transit orbits on the Poincar\'e section $\Sigma$, where the dots are the initial conditions for the corresponding trajectories. The right figure shows the transition tube for a given energy. The critical surface divides the transition tubes into two parts whose left part gives the initial conditions for orbits transitioning to the right, and right part gives the initial conditions for orbits transitioning to the left. Some trajectories are given to show how the transition tube controls the transition whose initial conditions are shown as dots on the left Poincar\'e section with same color.
	}
	\label{linear_tube_plot}
\end{figure}

\begin{figure}[!t]
	\begin{center}
		\includegraphics[width=\textwidth]{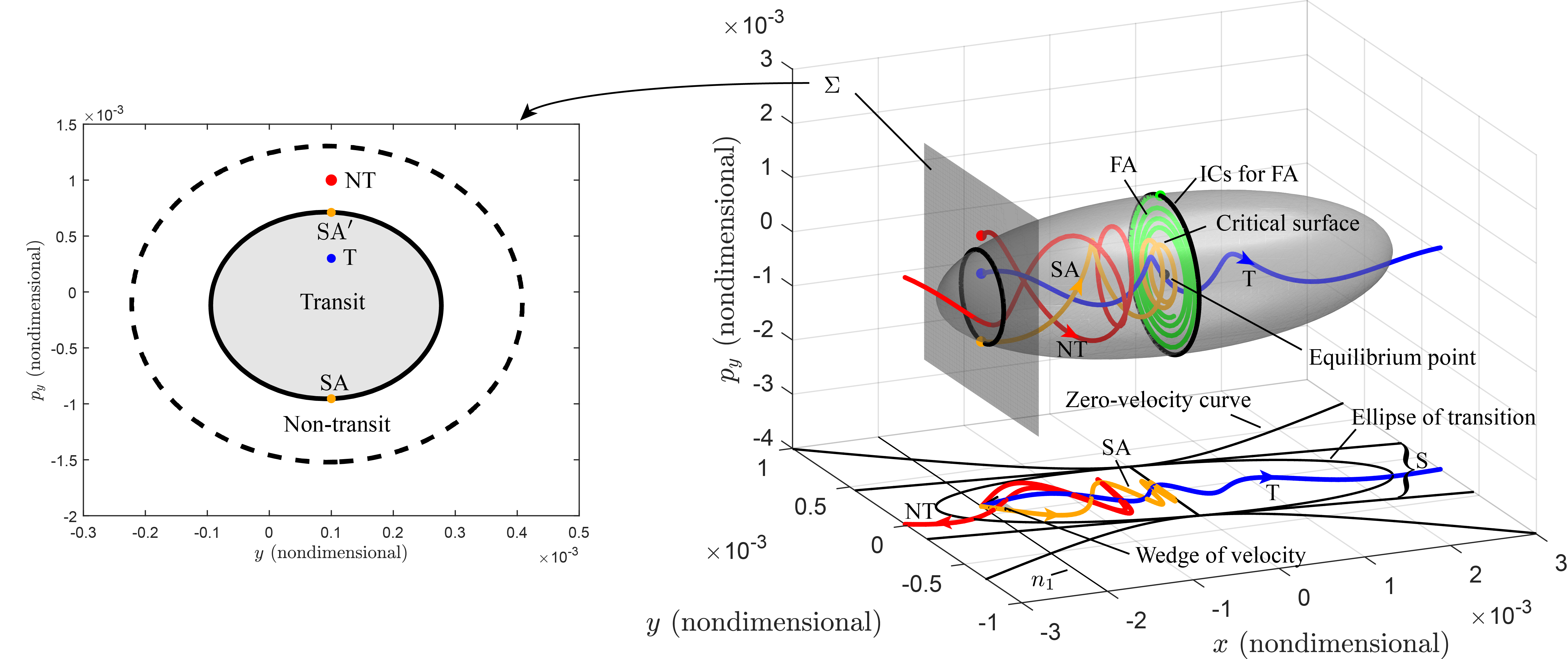}
	\end{center}
	\caption{\footnotesize 
		Transition ellipsoid for the dissipative system of initial energy $h$. The left figure shows the Poincar\'e section $\Sigma$, where the dots are the initial conditions for the corresponding trajectories with the same color in the right figure and the solid ellipse is the set of initial conditions for saddle-type asymptotic orbits. For comparison, the dashed ellipse of the tube boundary for the conservative system with the same energy $h$ is also given. On the right is the ellipsoid giving the initial conditions for all transit orbits. The critical surface divides the ellipsoid into two parts. Each side of the ellipsoid gives the initial conditions of transit orbits passing through the critical surface to the other side. In this figure, SA and FA denote the saddle-type and focus-type asymptotic orbits, respectively. Notice that due to the energy dissipation here, the periodic orbit in the conservative system becomes the initial conditions (ICs) of the focus-type asymptotic orbits.
	}
	\label{linear_ellipsoid_plot}
\end{figure}

\paragraph{Transition tube and transition ellipsoid.}
We have obtained the transition boundaries for both the conservative and dissipative systems. Their geometries in the phase space are in the form of a tube and an ellipsoid given in \eqref{tube_formula} and \eqref{ellipsoid formula} which are referred to as the transition tube and transition ellipsoid, respectively, in escape dynamics \cite{zhong2020geometry}. See the tube and ellipsoid in Figure \ref{linear_tube_plot} and Figure \ref{linear_ellipsoid_plot}, respectively.

In the figures, the transition tube and transition ellipsoid serve as the transition boundary in the phase space giving the initial conditions of a given energy $h$ for all the trajectories that can transit from the one side of the saddle point to the other in the conservative and dissipative systems, respectively. All the transit orbits must have the initial conditions inside the transition boundary, while the non-transit orbits have the initial conditions outside the boundary; surely, the transition boundary gives the initial conditions for the asymptotic orbits. To physically show the transition criteria, three initial conditions (one inside, one outside, and one on the transition boundary) are selected on the Poincar\'e section $\Sigma$ for both the conservative and dissipative systems. Three different types of trajectories (i.e., transit orbit, non-transit orbit, and asymptotic orbit) are observed. We can know that their transition conditions are truly governed by the transition boundary which demonstrates the rightness of the transition criteria and transition boundary we obtained.

Notice the transition tube and transition ellipsoid are divided into two parts by a surface, referred to as the \textit{critical surface}. The left part gives the initial conditions for the transit orbits going to the right and the right part gives the initial conditions for the transit orbits going to the left. The transit orbits can cross the critical surface, while the non-transit orbits will bounce back to the region where they come from. The analytical solution of the critical surface can show that the conservative system and dissipative system share the same critical surface. It means the critical surface is an intrinsic property of the system and will not be affected by any dissipation. In fact, the critical surface is a separatrix of different potential wells. Thus, once a transit orbit cross the critical surface, it can not return.

At the bottom of both Figure \ref{linear_tube_plot} and Figure \ref{linear_ellipsoid_plot}, shown is the flow on the position space (or configuration space) projected from the phase space. Notice that all possible motion is contained within the zero-velocity curve (corresponding to $\bar p_1=\bar p_2=0$), the boundary of  motion in the position space for a given value of energy. 
The projection of the transition tube and \textit{transition ellipsoid} onto the position space is a strip and an ellipse given in \eqref{strip} and \eqref{ellipse}, respectively. The strip is denoted by S. The ellipse in tube dynamics is called the \textit{ellipse of transition} \cite{zhong2018tube,zhong2020geometry}. Outside of the strip and ellipse, no escape is allowed. Even in the interior area, the escape is not definitely guaranteed. At each position inside of the strip and ellipse, there exists a wedge of velocity dividing the transit and non-transit orbits. The wedge of velocity gives the ``right" directions for the transition. Orbits with the velocity interior to the wedge are transit orbits, while orbits with the velocity outside of the wedge are non-transit orbit. Of course, orbits with velocity on the boundary are the asymptotic orbits. In fact, the boundaries of the wedge of velocity at a specific position are the lower and upper bounds of the velocity at that point. See A and A$'$ in Figure \ref{linear_tube_plot}, and SA and SA$'$ in Figure \ref{linear_ellipsoid_plot} for the lower and upper bounds of the velocity at a specific position in the conservative and dissipative systems, respectively. We can also observe that the periodic orbit in the position space is a straight line which means the periodic orbit in the phase space is perpendicular to the position space.

From Figure \ref{linear_tube_plot} and Figure \ref{linear_ellipsoid_plot}, we observe the transition tube encompasses the transition ellipsoid and they are tangential to each other at the critical surface. It means the dissipation in the system decreases the possibility of the transition. This is an intuitive conclusion.  In fact, the transit orbit in Figure \ref{linear_tube_plot} and the non-transit orbit in Figure \ref{linear_ellipsoid_plot} have the same initial condition. However, the dissipation makes a transit orbit in the conservative system become a non-transit orbit in the dissipative system. It demonstrates the dissipation in the system decreases the possibility of transition. When the dissipation is considered in the system, the transit orbit must start from a position not far from the equilibrium point; otherwise, the evanescent energy of the orbit will fall below the critical energy that allows the transition before crossing the critical surface. In this condition, the bottleneck around the index-1 saddle is closed so that the transition is impossible. The farthermost position for the transit orbit in the dissipative systems are the end points of the transition ellipsoid.


\paragraph{Stable global invariant manifold of the equilibrium point as separatrix.} 
In the discussion about the linearized dynamics around the index-1 saddle, the key observation  is that the initial conditions for the asymptotic orbits are on the surface of transition tube and transition ellipsoid in the conservative and dissipative systems, respectively. The tube and ellipsoid are constant Hamiltonian energy slices of points of the stable invariant manifolds of the periodic orbit and the equilibrium point, respectively; that is, the stable invariant manifolds can be foliated by the Hamiltonian energy $h$. For a given energy $h$, the stable invariant manifold separates two distinct types of motion: transit orbits and non-transit orbits. The transit orbits, passing from one realm to another, are those inside the invariant manifold. The non-transit orbits, which bounce back to their realm of origin, are those outside the invariant manifold. 

The concept of the invariant manifold of a periodic orbit and an equilibrium point will be important for the computation of transition boundary in the nonlinear system. In the next part, we will describe the algorithms to compute the invariant manifold as the transition boundary in the nonlinear system.

\section{Algorithms for computing the invariant manifolds}

\label{Manifold_algrithms}

\subsection{Invariant manifold of a periodic orbit}
\label{Manif_PO}
In this part, we aim to introduce the process to compute the invariant manifold of a periodic orbit. It has two separate parts. The first part concerns the algorithm for computing a periodic orbit, whereas the second concerns the computation of the stable and unstable manifolds of the periodic orbit.

\paragraph{Periodic orbits.}
A solution of the dynamical system \eqref{general ODEs} is a periodic orbit \cite{seydel2009practical} if there exists a least time interval $T>0$ which satisfies $x(t+T)=x(t)$ for all $t$. 
We will refer to the periodic trajectory as $\bar x(t)$. 
Multiple methods have been developed to compute the periodic orbits, such as the method of multiple scales \cite{nayfeh2008nonlinear,nayfeh2008applied}, incremental harmonic balance method \cite{lau1981amplitude,fu2006analysis}, and differential correction (or shooting method) \cite{KoLoMaRo2011,parker2012practical,sundararajan1997dynamics}, to name but a few. In the following, we will introduce another efficient BVP
approach which can compute the periodic orbits very accurately. Before discussing this approach, we rescale the time by introducing the linear transformation, $\tau=t/T$, so that the period $T$ appears explicitly in the equations of motion. Thus, the equations of motion in \eqref{general ODEs} can be rewritten by
\begin{equation}
\frac{\mathrm{d} x}{\mathrm{d} \tau}=T f(x), \hspace{0.3in} 0 \leqslant \tau \leqslant 1.
\label{rescaled ODEs}
\end{equation}
where $T$ is the unknown period. For the periodic orbits, we have the periodicity condition,
\begin{equation}
x(0)=x(1).
\label{periodicity condition}
\end{equation}
However, \eqref{rescaled ODEs} and \eqref{periodicity condition} do not uniquely determine the periodic solution, since if $x(t)$ is a periodic solution, so is $x(t + \delta)$. To avoid the arbitrary phase shift $\delta$, the following integral phase condition \cite{seydel2009practical,krauskopf2007numerical,dankowicz2013recipes} is widely used,
\begin{equation}
\int_0^T \left[x(t) - x^*(t) \right]^T f(x(t)) \mathrm{d}t=\int_0^1 \left[x(\tau) - x^*(\tau)\right]^T f(x(\tau)) \mathrm{d} \tau=0,
\label{PO_phase_condition}
\end{equation} 
where $x^*(t)$ is a known nearby solution.
The BVP 
is now formulated, 
and will require  numerical methods. 
The Matlab-based software package COCO \cite{dankowicz2013recipes} was applied to compute the periodic orbits.
COCO is a continuation tool which contains the algorithm described here as a toolbox, \textit{po}, which uses the collocation  and pseudo-arclength methods. 

\paragraph{Invariant manifold of a periodic orbit.} 
As mentioned in the introduction, the general way of computing the global invariant manifold is to globalize the local invariant manifold of the corresponding linearized system. Thus, here we can first find the local approximations of the manifold of the periodic orbit from the eigenvectors of the monodromy matrix and then grow the linear approximations by integrating the nonlinear equations of motion \eqref{general ODEs}. The procedure is known as globalization of the manifolds. Before growing the invariant manifold of the periodic orbit, we need to compute the state transition matrix $\Phi(t)$ along the periodic orbit which can be obtained by numerically solving the following variational equations from time $0$ to $T$,
\begin{equation}
\dot \Phi(t)= D f(\bar x(t)) \Phi (t), \hspace{0.2in} \text{with } \Phi(0)=I_n.
\end{equation}
Once the monodromy matrix $M\equiv\Phi(T)$ is obtained, its eigenvalues (the Floquet multipliers) can be computed numerically. For the two-mode Hamiltonian system about the shallow arch 
in the current study, $M$ is an infinitesimally symplectic matrix and its four eigenvalues 
consist of one real pair and one imaginary pair on the unit circle (see  \cite{KoLoMaRo2011}),
\begin{equation}
\lambda_1>1, \hspace{0.2 in} \lambda_2 =\tfrac{1}{\lambda_1}, \hspace{0.2 in} \lambda_3=\lambda_4=1.
\end{equation}

The eigenvector associated with eigenvalue $\lambda_1$ is in the unstable direction, while the eigenvector associated with eigenvalue $\lambda_2$ is in the stable direction. Denote the stable and unstable eigenvectors at the initial condition $x_0$ on the periodic solution by $e^s(x_0)$ and $e^u(x_0)$, respectively, normalized to unity. In this setting, we can obtain the initial guess for the stable and unstable manifolds, denoted by $x^s(x_0)$ and $x^u(x_0)$, at $x_0$ along the periodic orbit written in the following form,
\begin{equation}
\begin{aligned}
& x^s(x_0)=x_0 + \varepsilon e^s (x_0),\\
& x^u(x_0)=x_0 + \varepsilon e^u (x_0),
\end{aligned}
\end{equation} 
where $\varepsilon$ is a small parameter to obtain a small displacement from $x_0$ in the appropriate direction. 
The magnitude of $\varepsilon$ should be small enough so that the linear estimate can satisfy the accuracy, yet not so small that the time to obtain the global manifold becomes large due to the asymptotic behavior of the stable and unstable manifolds \cite{KoLoMaRo2011}. 

Once the initial guess for the stable and unstable manifolds at $x_0$ is obtained, it is straightforward to globalize the manifold. By numerically integrating the unstable vector forwards in time, using both $\varepsilon$ and $-\varepsilon$, one generates  trajectories shadowing the two branches, $W^{u+}$ and $W^{u-}$, of the unstable manifold of the periodic orbit. Similarly, by integrating the stable vector backwards, we generate a trajectory shadowing the two branches of the stable manifold, $W^{s \pm}$. 
For a trajectory on the manifold at some other point $\bar x(t)$ on the periodic orbit, one can  use the state transition matrix to transport the eigenvectors from $x_0$ to $\bar x(t)$,
\begin{equation}
e^s(\bar x(t))= \Phi(t) e^s(x_0), \quad \text{and} \quad e^u(\bar x(t))= \Phi(t) e^u(x_0).
\end{equation}
Since the state transition matrix does not preserve the norm, the resulting vectors must be renormalized. 

Globalizing the manifold at $N$ points  $\bar x(t)$ on the periodic orbit (where $N$ is large) provides a set of $N$ trajectories which approximate the global manifold of energy $h$. 
In the case of the stable manifold, one therefore obtains the boundary of the transit orbits starting with energy $h$, $\partial \mathcal{T}_h$.

\subsection{Invariant manifold of an equilibrium point}

In the previous section, we discussed the approach to compute a periodic orbit of energy $h$ and its stable and unstable invariant manifolds in the conservative system. The stable invariant manifold along each energy manifold of energy $h$ is $\partial \mathcal{T}_h$, the boundary of the initial conditions starting at energy $h$ that will soon escape from one side of the index-1 saddle to the other. 
Once this is understood, it is 
natural to consider
what the global phase space structure 
governing the transition will be in the dissipative system. 
We address that concern by computing the invariant manifold of the equilibrium point in the dissipative system.

We consider the same general form of a dynamical system in \eqref{general ODEs} to define the dissipative system with a hyperbolic equilibrium point, 
$x_e$. 
The Jacobian of the equilibrium point, $Df(x_e)$, has $k$ eigenvalues with negative real part. 
The real parts of the $k$ eigenvalues and the corresponding generalized eigenvectors are written by $\lambda_i^s<0$ and $u_i$ $(i=1,\cdots,k)$, respectively. 
Thus, the saddle has a $k$-dimensional local, invariant stable manifold, denoted by $W_{loc}^s(x_e)$, which is tangent to the respective invariant stable subspaces, $E^s$, of the linearized system about the saddle, spanned by the stable eigenvectors $u_i$. 

Once the local stable manifold is determined, it can be globalized to obtain the global $k$-dimensional stable manifold $W^s (x_e)$ \cite{krauskopf2006survey}. 
The direct approach to obtain the global manifold is to select initial conditions in the stable subspace  a small distance from the equilibrium point 
and integrate backward in time, thereby obtaining  orbit segments on the stable manifold. 
Numerical continuation by using the  resulting orbit as a starting solution might give the global manifold. 
However, some challenges may appear  \cite{krauskopf2003computing}, such as  large aspect ratios of the computed manifold surface due to a difference in the real parts of the eigenvalues, and corresponding stretching of the distance between solutions after a sufficiently long integration. 
To solve these problems, re-meshing of the  manifold surface is needed, which is another challenge.

Another way of obtaining the stable manifold is solving a proper two-point BVP 
\cite{krauskopf2007numerical} which can control the endpoints of the trajectories. Before describing the process, we need to rescale the time by $t=T \tau$ which puts \eqref{general ODEs} into the same form in \eqref{rescaled ODEs} where $\tau$ varies from 0 to 1. It should be mentioned that, compared to the period $T$ of a periodic orbit in the conservative system, $T$  in the dissipative case is a chosen  time-scale. 
We can consider $T$ as either a parameter or a function whose derivative with respective to $\tau$ is zero, i.e., $d T /d \tau=0$. Here we will use the latter one. To form a complete BVP, 
we still need some boundary conditions. 
The boundary conditions at $\tau =0$ can be selected on an initial hyper-sphere on the stable subspace given by,
\begin{equation}
{x}(0)={x}_e + r_0 \sum_{i=1}^{k} a_i {u}_i,
\label{general initials}
\end{equation}
where $a_i$ are  parameters controlling the direction of the initial condition; 
$r_0$ is the distance of the initial conditions from the equilibrium point. 
The parameter $r_0$, like $\varepsilon$ in the previous section,  should be properly selected, neither too small nor too large.

In the following, we will take the snap-through of a shallow arch with damping as an example. The equations for the 
BVP
to compute the invariant manifold of the index-1 saddle are,
\begin{equation}
\begin{aligned}
\dot X &=T\frac{p_X}{M_1},  \hspace{0.2in} && \dot Y = T\frac{p_y}{M_2}, \\
\dot p_X &= T \left(- \frac{\partial \mathcal{V}}{\partial X} - C_H p_X \right), \hspace{0.2in}
&&\dot p_Y = T \left( - \frac{\partial \mathcal{V}}{\partial Y} - C_H p_Y \right),  \\
\dot T&=0.
\label{eom_manifold_damp}
\end{aligned}
\end{equation}
where $\partial \mathcal{V}/ \partial X$ and $\partial \mathcal{V}/ \partial Y$ are given by in \eqref{eq:eomHam}.

In the dissipative system, the index-1 saddle has become a hyperbolic point with a $k=3$-dimensional stable invariant manifold. 
The $\tau = 0$ boundary conditions 
can be selected along an initial  2-sphere with radius $r_0$ given by,
\begin{equation}
{x}(0)={x}_e + r_0 \left(\sin \theta \sin \phi u_1 + \sin \theta \cos \phi {u}_2 + \cos \theta {u}_3 \right)
\label{initials_for_arch}
\end{equation}
where $\theta$ and $\phi$ are the two parameters (spherical coordinates). 
Notice that, for fixed $r_0$ and for each $\theta$ $\phi$, \eqref{initials_for_arch} corresponds to 4 boundary conditions at $\tau=0$. 
We still need one more boundary condition, at $\tau =1$. 
We can have several choices, such as the energy, arclength or time of the trajectories. 
In the current problem we want to find the  boundary in the dissipative system of transition trajectories with an initial Hamiltonian energy $\mathcal{H}=h$. 
This is done by assigning the energy to the endpoint at $\tau=1$, 
\begin{equation}
\mathcal{H} \left( X(1),Y(1),p_{X}(1),p_{Y}(1) \right) =h,
\label{manifold-energy}
\end{equation}
After we set up the BVP, 
we can apply numerical continuation to obtain the invariant manifold of the system. 
Before solving the BVP, 
we need to prescribe $r_0$ and $h$, so we have 5 variables ($X$, $Y$, $p_X$, $p_Y$, and $T$) and 2 parameters ($\theta$ and $\phi$). 
On the other hand, we have 5 boundary conditions in \eqref{initials_for_arch} and \eqref{manifold-energy}. Thus, the BVP here is a two-parameter continuation. To simplify the continuation process, we can reduce the system to a one-parameter continuation by introducing a proper Poincar\'e section.

Of the three stable eigenvectors ${u}_1$, ${u}_2$, and ${u}_3$, we assume the magnitude of the real part of the eigenvalue associated with ${u}_1$ is the largest. 
Thus, the   ${u}_1$ direction is the dominant stable direction. 
If we take $\theta=\phi =\pi/2$ in \eqref{initials_for_arch} and use it as an initial condition to numerically integrate \eqref{eom_manifold_damp} {\it backward in time} until the trajectory reaches the desired energy $h$, we  reach the furthest end of the nonlinear transition ellipsoid. 
Since this trajectory is approaching the saddle in positive time along the most stable direction, it is the fastest stable asymptotic orbit to the saddle. We refer to it as the \textit{fastest trajectory}. 
In the following we will use it as a reference trajectory to determine an appropriate Poincar\'e section. 

Let $\mathbf{r}$ denote the position vector of an arbitrary point on the fastest trajectory in the $X$-$Y$-$P_Y$ subset of phase space. We can obtain the tangent vector $\mathbf{t}$ at that point along the fastest trajectory,
\begin{equation}
\mathbf{t}=\frac{\partial \mathbf{r}}{\partial s}= \frac{\partial X}{\partial s} \mathbf{e}_X + \frac{\partial Y}{\partial s} \mathbf{e}_Y + \frac{\partial p_Y}{\partial s} \mathbf{e}_{p_Y}= \frac{\dot X}{\dot s} \mathbf{e}_X +\frac{\dot Y}{\dot s} \mathbf{e}_Y +\frac{\dot p_Y}{\dot s} \mathbf{e}_{p_Y},
\end{equation}
where $\mathbf{e}_X$, $\mathbf{e}_Y$, and $\mathbf{e}_{p_Y}$ are the corresponding basis vectors long $X$, $Y$, and $p_Y$. Here $s$ is the arclength of the reference trajectory in $X$-$Y$-$P_Y$ space, which is a function of $\tau$, defined by,
\begin{equation}
s(\tau) =\int_0^\tau \sqrt{\dot X^2 + \dot Y^2 + \dot p_Y^2} \mathrm{d} \tau^{\prime}. 
\end{equation}

\begin{figure}[!h]
	\begin{center}
		\includegraphics[width=\textwidth]{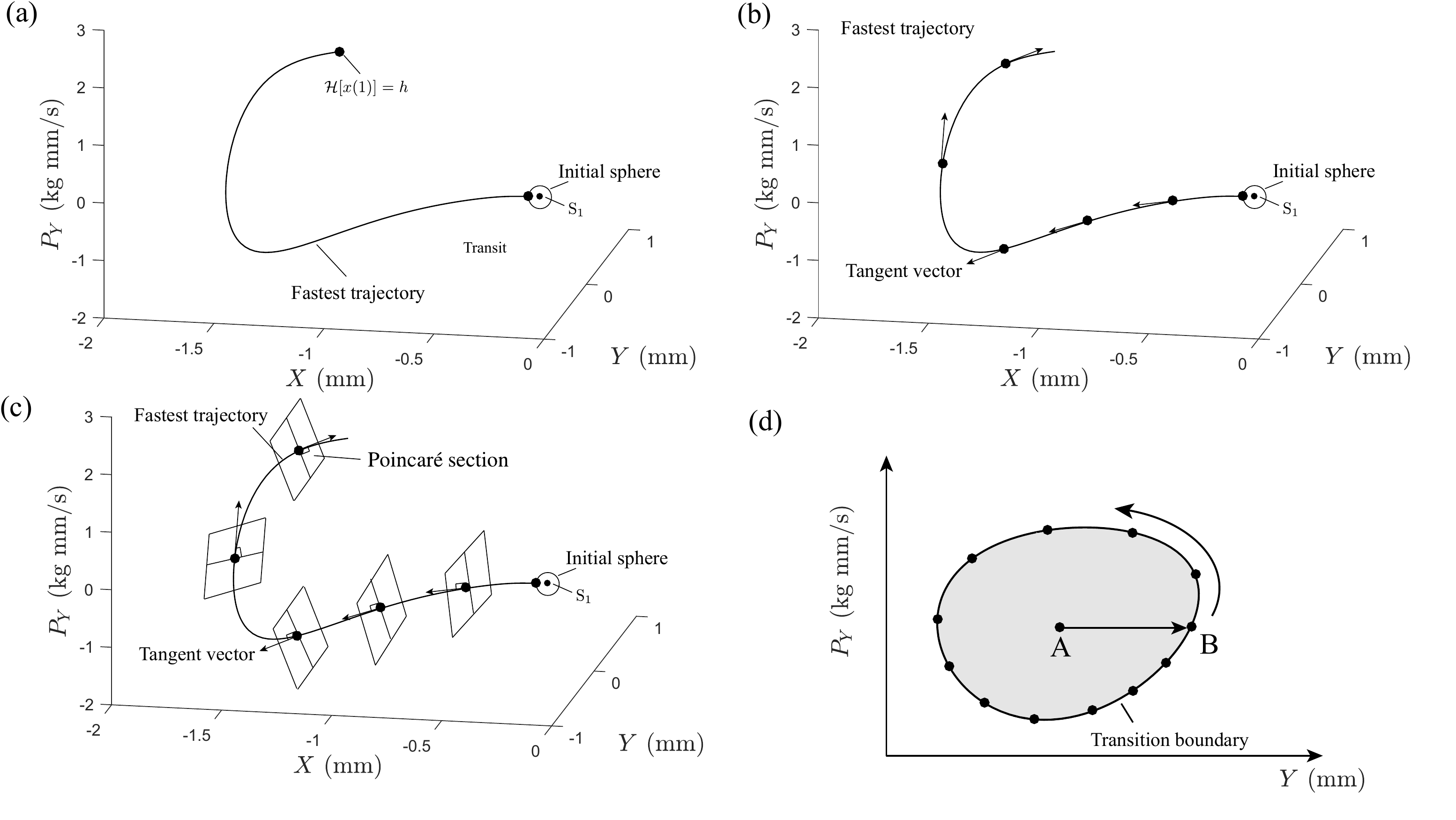}
	\end{center}
	\caption{\footnotesize 
		Illustration of selecting an extra Poincar\'e section to reduce the two-parameter continuation to a one-parameter continuation: 
		(a) Select the initial condition of the fastest trajectory on the initial sphere (with small radius $r_0$) in the stable subspace of the linearized system and numerically integrate the nonlinear equations until the trajectory reaches the given energy $h$. 
		(b) Select a bunch of points on the fastest trajectory and compute the tangent vector at each point along the fastest trajectory. Each point has uniform arc-length to its two neighboring points; 
		(c) Finally the plane normal to the tangent vector at each point can be selected as the Poincar\'e section at that point. 
		(d) After we determine the Poincar\'e sections, we can select another initial condition on the initial sphere and numerically integrate the nonlinear equations until it reaches the Poincar\'e sections, denote the intersection as point A. Of course, we can also use the fastest trajectory. In general, point A lower than the given energy $h$. Next, we can fix $p_Y$ and commit the continuation along $Y$ direction until the Hamiltonian reaches $h$ so that we can obtain the point B which is on the transition boundary. Then we can use point B as the starting solution and do the continuation with fixed total energy $h$ by which we can obtain the transition boundary on the Poincar\'e section.
	}
	\label{Selection_POS}
\end{figure}

For a specific point $(X_0,Y_0,p_{Y_0})$ on the fastest trajectory, we can choose a plane normal to $\mathbf{t}$ at that point as the Poincar\'e section. The mathematical expression of the plane is given by,
\begin{equation}
t_X (X(1)-X_0) + t_Y (Y(1)-Y_0) + t_{p_Y} (p_Y(1)-p_{Y_0})=0,
\end{equation}
where $t_X$, $t_Y$, and $t_{p_Y}$ are the components of the tangent vector 
$\mathbf{t}$ along the $X$, $Y$, and $p_Y$ axes, respectively. In this case, we have one more algebraic equation as the extra boundary condition. 
This means our problem reduces to one-parameter continuation on the Poincar\'e section. 
Figure \ref{Selection_POS} gives the illustration of the process to select the proper Poincar\'e section described above. The algorithm of solving the boundary-value problem will be implemented in COCO \cite{dankowicz2013recipes} again to compute the invariant manifold of the equilibrium point. 
In this way one  obtains the transition boundary, $\partial \mathcal{T}_h$, in the dissipative system.

\section{Numerical results}
In this section, we give the geometry of transition boundary that mediates the nonlinear snap-through buckling of a shallow arch in both the conservative and dissipative systems. In the corresponding examples, the geometrical and material parameters are selected following previous   experimental and theoretical studies \cite{WiVi2016,zhong2018tube}: 
$b=12.7$ mm $d=0.787$ mm, $L=228.6$ mm, $\gamma_1 = 0.082$ mm, and $\gamma_2 = -0.077$ mm; 
the Young's modulus and the mass density are $E=153.4$ GPa and $\rho=7567 \ \mathrm{kg \ m^{-3}}$; moreover, the thermal load is chosen to be $184.1$ N. 
For the convenience of discussing the energy, we use the {\it excess energy} $\Delta E$ \cite{naik2019finding} above the saddle point S$_1$ which is defined by $\Delta E=E - E_c$. The  energy of S$_1$ is $E_c$, the critical (minimum) energy necessary for  transition between the two wells. 
For positive excess energy, $\Delta E>0$, the bottleneck region around the saddle is open so that the trajectories have a chance to escape; otherwise the bottleneck region is closed and  transition is not energetically possible.

\subsection{Conservative systems}

\begin{figure}[!t]
	\begin{center}
		\includegraphics[width=\textwidth]{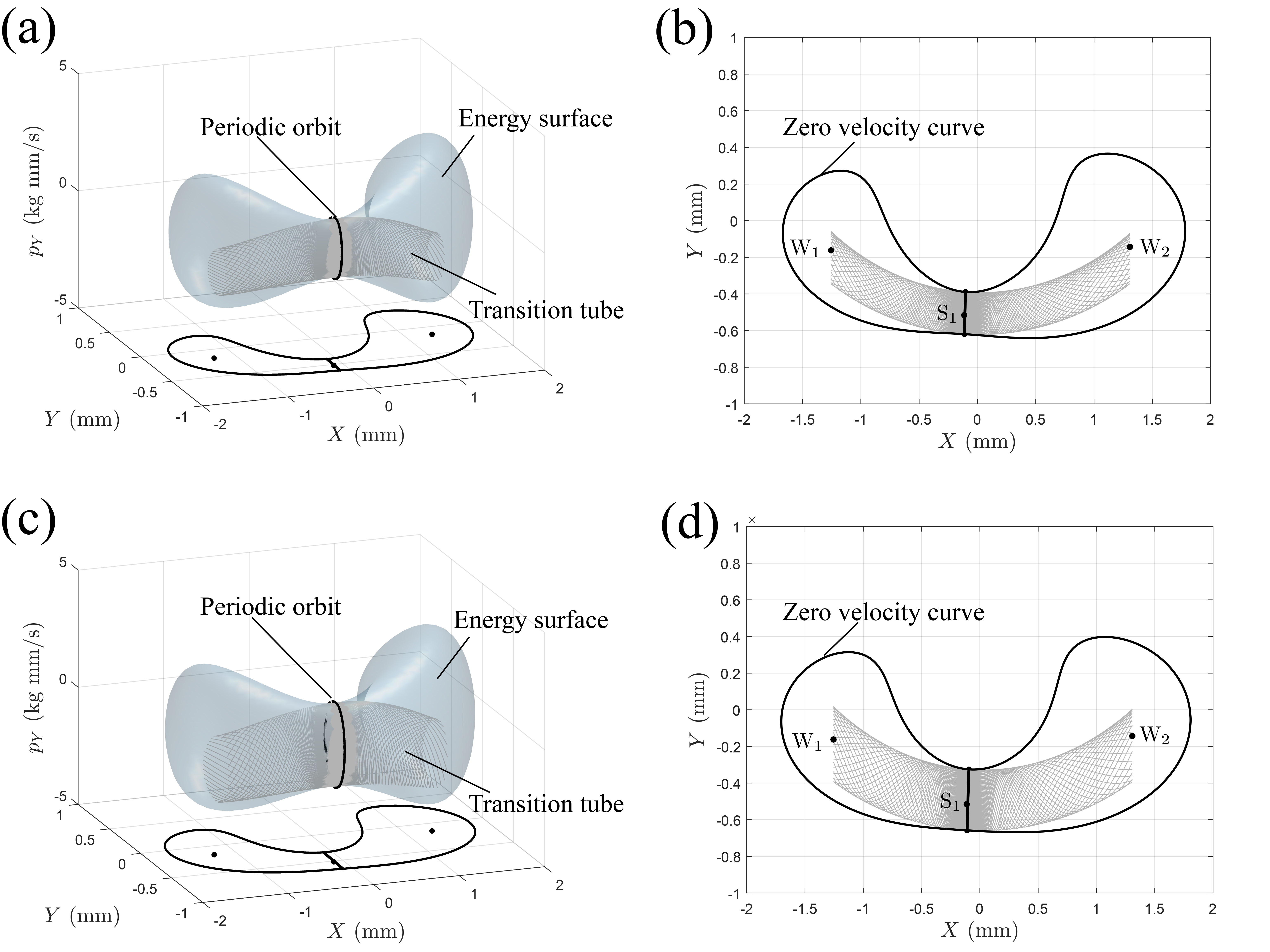}
	\end{center}
	\caption{\footnotesize 
		Transition tubes obtained by current algorithm in the conservative system. The left two in (a) and (c) are the transition tubes in the phase space with excess energy $\Delta E=1.0 \times 10^{-4}$ J and $\Delta E=2.0 \times 10^{-4}$ J, respectively. The right two in (b) and (d) are the corresponding projections onto configuration space. 
		Exterior to the transition tubes, the energy manifolds are shown which bound the possible motion of all trajectories with that initial energy. 
	}
	\label{multiple_nonl_tube}
\end{figure}

In this section, we give some examples of the transition boundary in the conservative system. 
Analogous to the linearized dynamics around S$_1$ as discussed in 
Section \ref{linear_conservative_system}, for initial conditions of a certain energy above the $E_c$, the transition between the potential wells 
in the 
conservative system is governed by a cylindrical conduit. 

Figure \ref{multiple_nonl_tube} shows two transition tubes with initial excess energy $\Delta E=1.0 \times 10^{-4} $ J and $\Delta E=2.0 \times 10^{-4}$ J, respectively. The left two are the transition tubes in the phase space and the right two are the corresponding projections onto the configuration space. Outside the transition tubes, we also plot the energy manifolds which is the boundary of all possible motions of the corresponding energy in the phase space. The transition tubes, i.e., the stable invariant manifolds of periodic orbits about the index-1 saddle point, are cylindrical tubes of trajectories asymptotically approaching  the periodic orbit in forward time.
The transition tube is the boundary in  phase space separating  transition and non-transition trajectories. All the trajectories of initial excess energy $\Delta E$ transitioning to a different potential well are inside  the tube manifold. 
The size of the transition tube compared to that that of the energy manifold is a measure of the probability of transition. Notice that this ratio is  larger for the case of larger excess energy. 
Moreover, since the energy in the conservative system keeps constant during evolution of any trajectories, the symplectic cross-section of the tube manifold is invariant, obeying Hamilton's canonical equations (with no dissipation). 
Notice that in the linearized system around the index-1 saddle, the transition tube 
appears as a straight cylinder. However, due to the nonlinear terms, 
the transition tube in the full system is curved.

\begin{figure}[!h]
	\begin{center}
		\includegraphics[width=\textwidth]{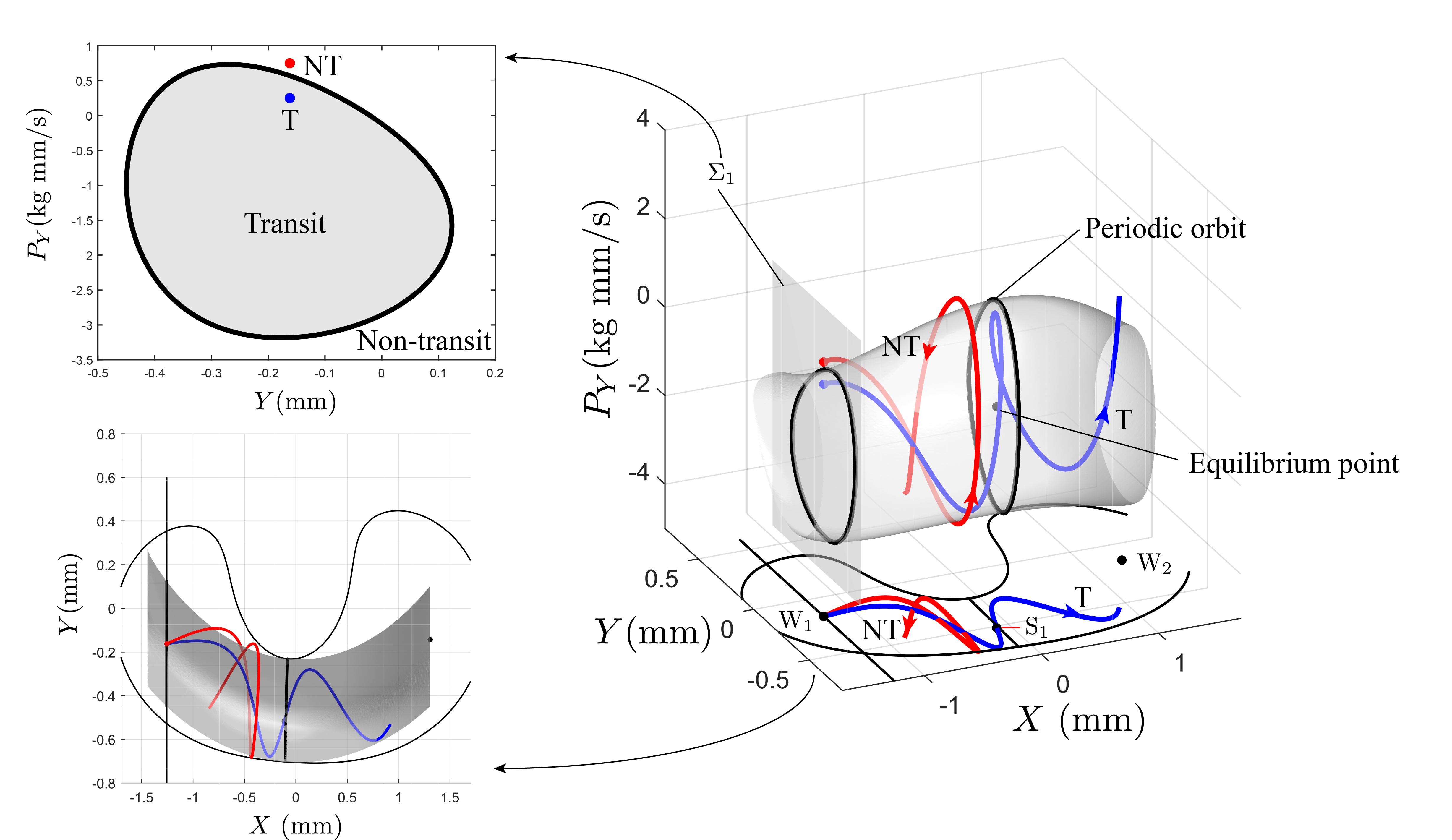}
	\end{center}
	\caption{\footnotesize 
		%
		A transition tube in the conservative system obtained by the boundary value problem approach: the right figure shows the transition tube in a 3-dimensional projection of the 4-dimensional phase space; the lower left shows the configuration space projection; the upper left shows the transition boundary, a closed curve, on the Poincar\'e section $\Sigma_1$ which separates the initial conditions with a given fixed energy for the transit and non-transit trajectories. 
		A transit orbit and a non-transit trajectory starting with initial conditions labeled by T and NT are shown, which are inside and outside of the transition boundary on the Poincar\'e section $\Sigma_1$, respectively.
		%
	}
	\label{Manifold_tube}
\end{figure}

To show how the transition tube confines the transition between potential wells, Figure \ref{Manifold_tube} gives the transition tube and two trajectories with $\Delta E=3.68 \times 10^{-4}$ J, which is coincident with the energy used in \cite{zhong2018tube}. 
A Poincar\'e section $\Sigma_1$ is selected which is defined by the $X$ value equal to that of the stable equilibrium point W$_1$. The intersection of the transition tube with $\Sigma_1$ 
is a closed curve serving as the transition boundary on the Poincar\'e section. 
At the location of W$_1$, two initial conditions are chosen, one inside and one outside of the transition boundary. The trajectory with the initial condition inside of the transition boundary is a transit orbit, transitioning from potential well W$_1$ to potential well W$_2$. The trajectory with the initial condition outside of the transition boundary, however, is a non-transit orbit, returning to W$_1$ before entering the realm of potential well W$_2$. 

\begin{figure}[!t]
	\begin{center}
		\includegraphics[width=0.5\textwidth]{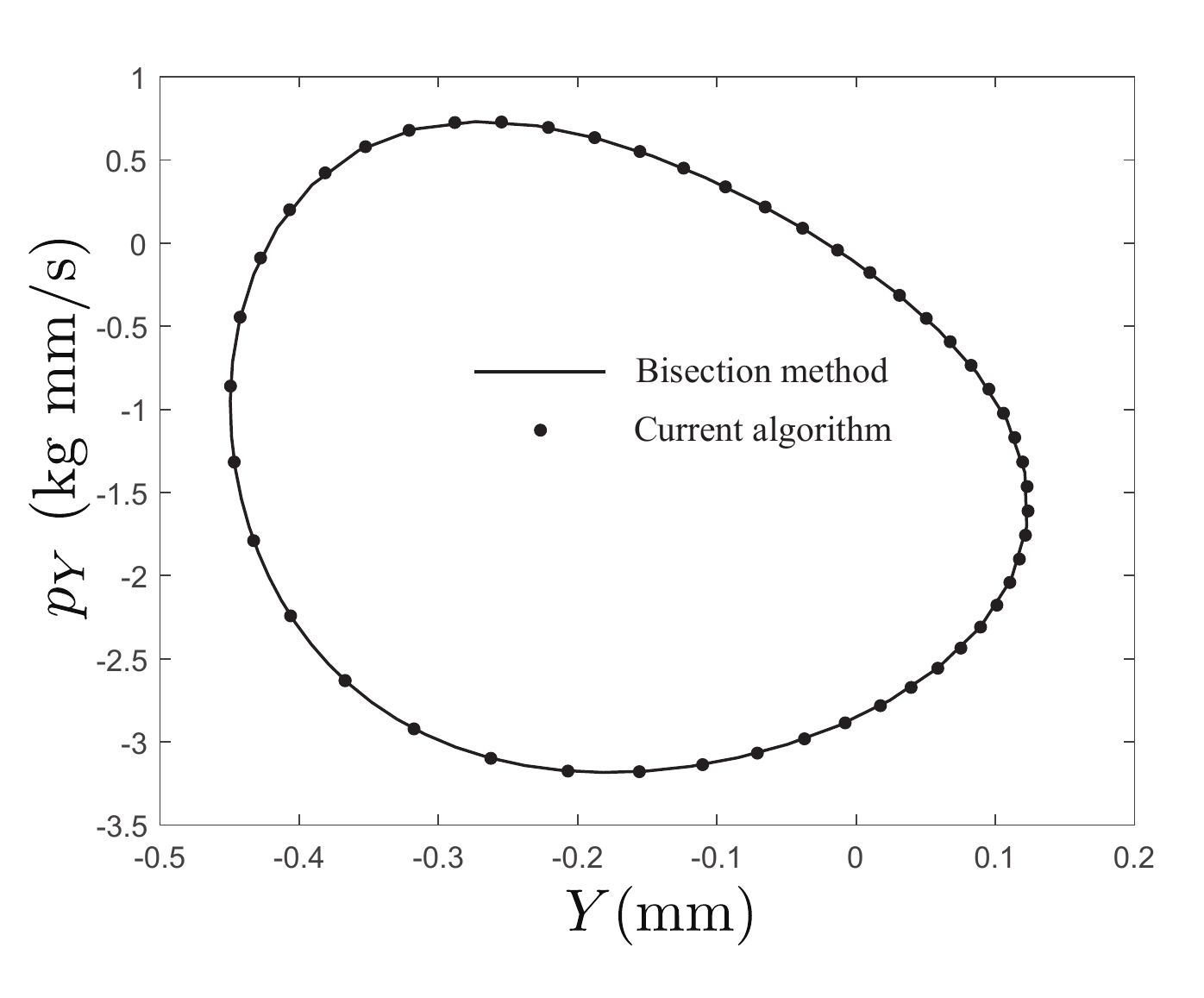}
	\end{center}
	\caption{\footnotesize 
		Comparison of the transition boundary on Poincar\'e section $\Sigma_1$ in the conservative system between the results obtained by the current algorithm and those obtained by the bisection method \cite{zhong2018tube}. The excess energy $\Delta E$ is selected to be $3.68 \times 10^{-4}$ J. The current result calculated via the boundary value problem approach is shown as dots and the result calculated by the bisection method \cite{zhong2018tube} is shown as a solid curve.
	}
	\label{conservative_comp}
\end{figure}

To further validate the current method, Figure \ref{conservative_comp} shows a comparison of the transition boundary on the Poincar\'e section $\Sigma_1$ between the current method and a previously developed bisection method \cite{zhong2018tube}. 
From the figure, the results obtained by the two methods agree well with each other, 
demonstrating that the two algorithms are consistent. 

The bisection method has an advantage over the BVP method in that it does not depend on information regarding the linearized dynamics about the transition. 
Instead, it is a `brute-force' approach which directly searches the boundary on a Poincar\'e section, 
and therefore has the disadvantage of taking a larger computational time for the same level of precision. 
A hybrid approach could be possible wherein one firstly obtains a boundary point along a specific direction on the Poincar\'e section using the bisection method. This point is an 
initial condition for a trajectory asymptotic  to the periodic orbit, and can be used to obtain the periodic orbit  itself. 
Once the periodic orbit is obtained, one can apply the globalization of the local invariant manifold of a periodic orbit via the BVP approach.


\subsection{Dissipative systems}

\begin{figure}[!t]
	\begin{center}
		\includegraphics[width=\textwidth]{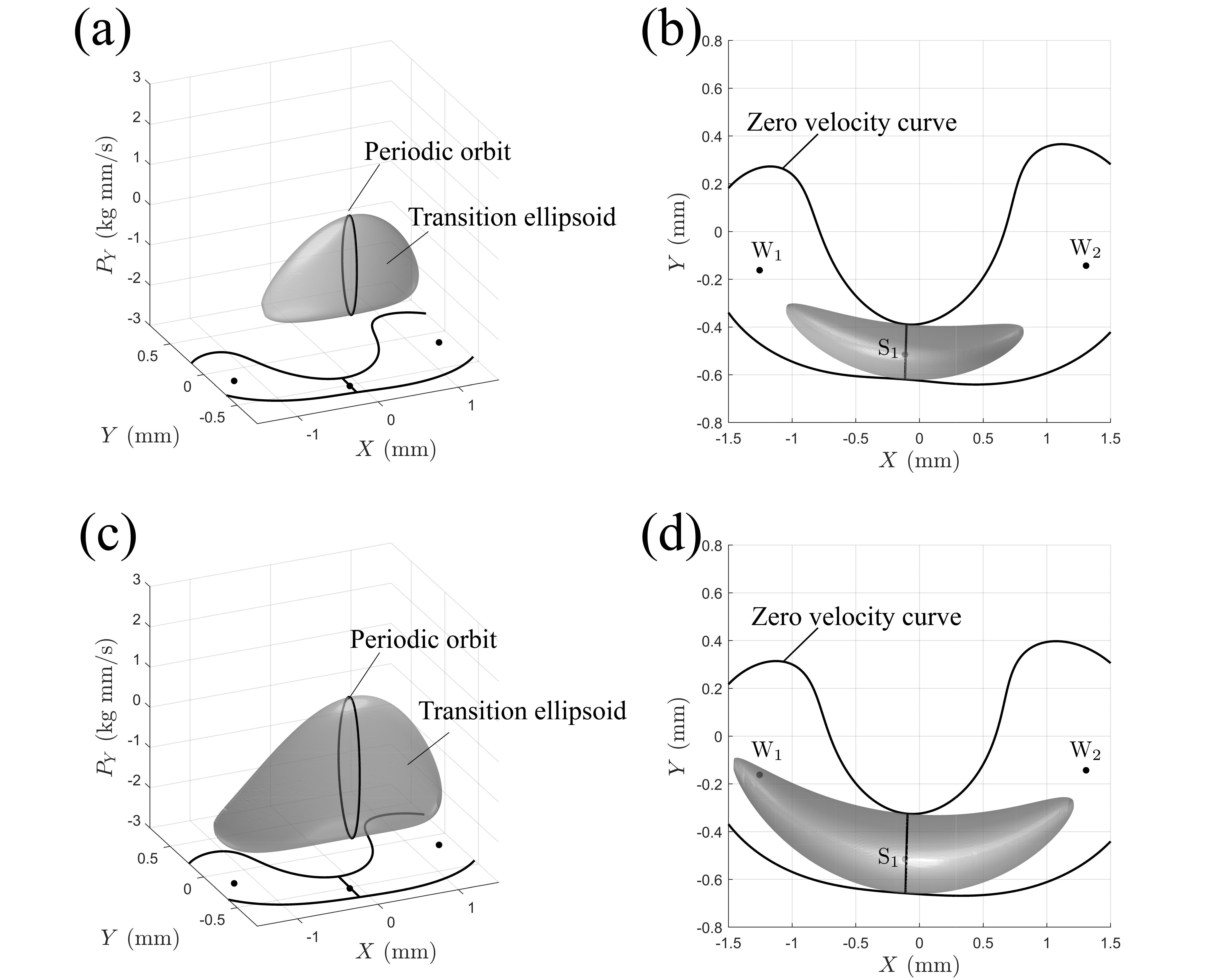}
	\end{center}
	\caption{\footnotesize 
		%
		Transition ellipsoids in the dissipative system obtained by the boundary value problem approach: (a) and (c) show the three-dimensional transition ellipsoids with excess energy $\Delta E = 1.0 \times 10^{-4}$ J and $\Delta E = 2.0 \times 10^{-4}$ J, respectively; (d) and (d) show the corresponding configuration space projections. An animation for the transition ellipsoid is at \url{https://www.youtube.com/watch?v=qzKQWe__uv4}
	}
	\label{Manifold_ellipsoid_1_2}
\end{figure}

In this section, 
the transition boundary in the dissipative system for the snap-through buckling of the shallow arch
is obtained. 
In the conservative system, the energy remains constant in time for all motions. 
However, in the dissipative system the energy  decreases as trajectories go forward in time. 
Furthermore, as in the linearized system, the phase space structure that governs the transition in the dissipative system is topologically distinct from that in  the conservative system. 
For the following numerical results in the dissipative system, 
the damping parameter is taken as $C_H=80$ s$^{-1}$.

Figure \ref{Manifold_ellipsoid_1_2} shows 
two transition ellipsoids with initial excess energy $\Delta E=1.0 \times 10^{-4} $ J and $\Delta E=1.0 \times 10^{-4}$ J, respectively. 
The corresponding configuration space projections are given on the right. 
The figure shows that the  transition ellipsoid of larger energy has a larger size 
relative to the energy manifold, compared with the smaller energy transition ellipsoid. 
That is, the probability for transition increases with initial excess energy. 
Due to the presence of nonlinear terms in the system, the transition ellipsoids appear curved compared to  their  
linearized system counterparts. 
However, their topology is the same: a 2-sphere. 

In Figure \ref{Manifold_ellipsoid_1_2}, the periodic orbits with the same excess energy from the 
conservative system are also shown. 
Notice that the periodic orbits are exactly on the boundary of the transition ellipsoids. 
In fact, the points on the periodic orbits are the initial conditions of the focus-type asymptotic orbits.
Each periodic orbit divides the corresponding transition ellipsoid into two parts. 
The left part of the transition ellipsoid bounds the initial conditions for transit orbits moving from  left well to the right well, while the right part bounds the initial conditions for  transit orbits moving from  right-to-left.

\begin{figure}[!t]
	\begin{center}
		\includegraphics[width=\textwidth]{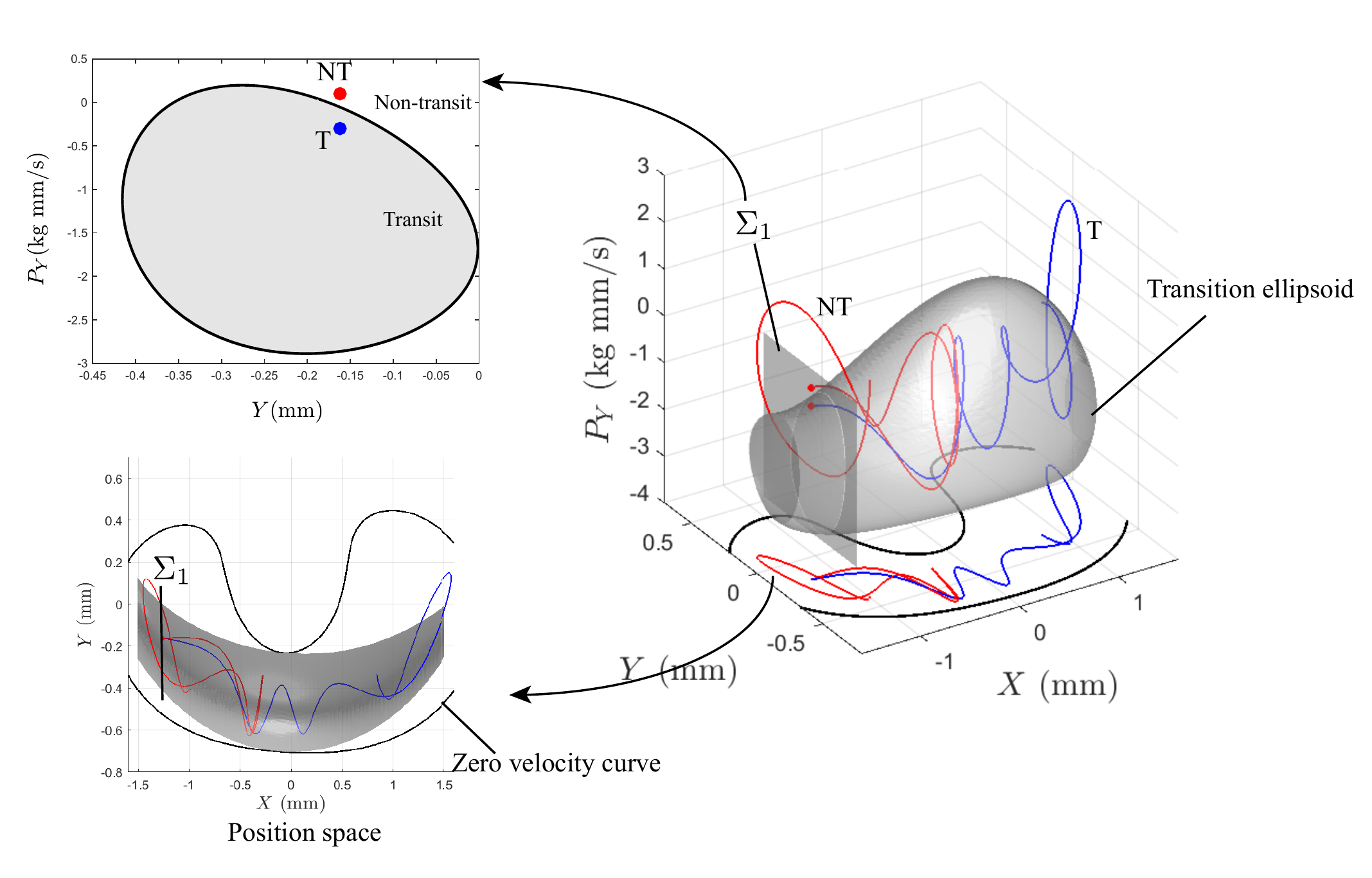}
	\end{center}
	\caption{\footnotesize 
		A transition ellipsoid in the dissipative system obtained by the boundary value problm approach: the right figure shows the transition ellipsoid in a 3-dimensional projection of the 4-dimensional phase space; the lower left shows the configuration space projection; the upper left shows the transition boundary, a closed curve, on the Poincar\'e section $\Sigma_1$ which separates the initial conditions with a given fixed energy for the transit and non-transit trajectories. A transit orbit and a non-transit trajectory starting with initial conditions labeled by NT and T are shown, which are inside and outside of the transition boundary on the Poincar\'e section $\Sigma_1$, respectively.
	}
	\label{Manifold_ellipsoid_368}
\end{figure}

\begin{figure}[!t]
	\begin{center}
		\includegraphics[width=0.5\textwidth]{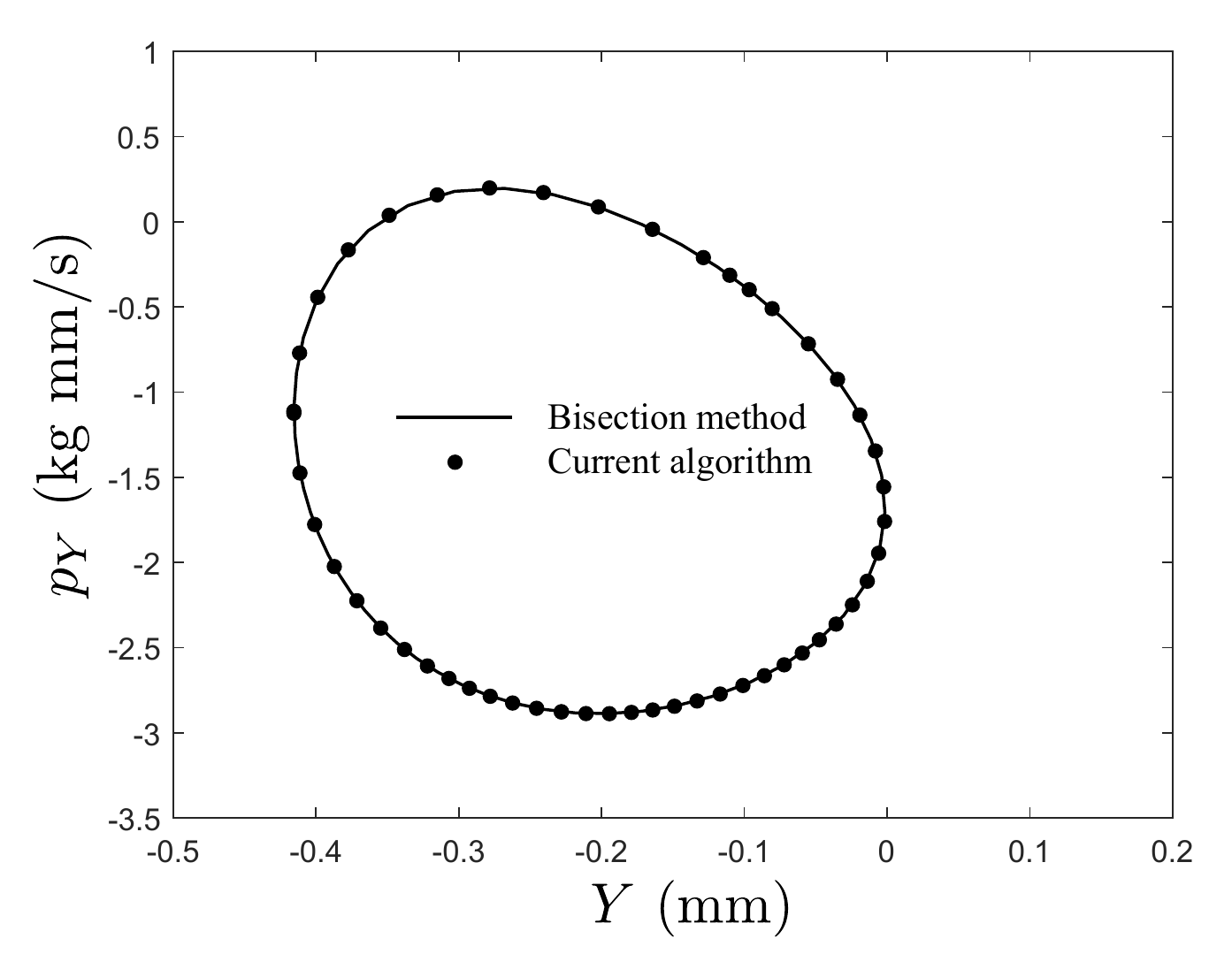}
	\end{center}
	\caption{\footnotesize 
		Comparison of the transition boundary on the Poincar\'e section in the dissipative system $\Sigma_1$ between the current algorithm and the bisection method \cite{zhong2018tube}. The excess energy is selected as $\Delta E=3.68 \times 10^{-4}$ J, and the coefficient of the linear viscous damping is taken as $C_H=80$ s$^{-1}$. The results obtained by current algorithm and the bisection method \cite{zhong2018tube} are shown by dots and solid curve, respectively. 
	}
	\label{dissipative_comp}
\end{figure}

To illustrate in further detail how the transition ellipsoid bounds 
the initial conditions leading to  transition, 
we select an excess energy $\Delta E=3.68 \times 10^{-4}$ and 
compute the transition ellipsoid  shown in
Figure \ref{Manifold_ellipsoid_368}. 
We select two initial conditions  on the Poincar\'e section $\Sigma_1$, both with a configuration space value equal to the equilibrium point W$_1$, but with non-zero velocity.
One initial condition is inside and the other outside of the transition ellipsoid boundary. 
Integrating the initial conditions forward in time, we  obtain two trajectories. 
From the figure, we find that trajectory T with the initial condition inside of the transition boundary escapes from the potential well W$_1$ to the other potential well W$_2$, 
while trajectory NT with the initial condition outside of the transition boundary bounces back to the region of origin. 
The same case study was conducted in \cite{zhong2018tube} via the bisection method. 
The comparison of the transition boundary on the Poincar\'e section $\Sigma_1$ obtained by the current study and \cite{zhong2018tube} is given in Figure \ref{dissipative_comp}. 
Good agreement between the two methods is observed. 


\section{Conclusion}

In this paper, we apply the concept of invariant manifolds to identify the boundaries of transition orbits in
a two degree of freedom nonlinear system with and without energy dissipation. 
The example system considered is 
the snap-through buckling of a shallow arch, where energy dissipation is necessary to model the behavior of the real system \cite{WiVi2016}. 
The essence of the snap-through buckling is the escape or transition from one potential well to another. 
The phase space structures that govern the transition in the conservative and dissipative systems are the stable invariant manifold of a periodic orbit and of an equilibrium point, respectively, of a prescribed energy. 
The global stable invariant manifolds are computed numerically by solving proper boundary-value
problems which are implemented in the continuation numerical package COCO \cite{dankowicz2013recipes}.

In the conservative system, the computational process providing the invariant manifold has two steps, first the computation of the periodic orbit by solving a proper boundary-value problem and second, the globalization of the manifold. 
In the dissipative system, we compute the invariant manifold of the index-1 saddle by another set of boundary-value problem: we first compute the stable subspace of the linearized system and 
select a hyper-sphere with small radius in this subspace.
The boundary conditions are selected as the points on the small sphere near the equilibrium point and the prescribed energy at the other end.
The boundary-value problem set-up for the dissipative system is a two-parameter continuation. 
To reduce the difficulty of conducting this continuation, 
a Poincar\'e surface-of-section is introduced so that the problem becomes a one-parameter continuation.

By using these approaches, one obtains the transition tube and transition ellipsoid serving as the transition boundary for the conservative and dissipative systems, respectively, 
which are topologically the same as those in the linearized dynamics. 
Trajectories with initial conditions inside the transition boundary will snap-through, while trajectories with initial conditions outside the transition boundary will not. 
As a demonstration for the accuracy and efficiency of the current algorithms, we compared the current results with those obtained by a bisection method used in \cite{zhong2018tube}. 
There is good agreement between each method, but the boundary value problem approach is more systematic. 

In \cite{zhong2020geometry}, the linearized dynamics underlying escape and transition in several widely known physical systems in the presence of dissipative and/or gyroscopic forces was summarized. 
The current study extends the linearized dynamics to the nonlinear case presenting the 
boundary-value problem approach to compute the transition boundary 
when dissipation is considered. 

Given the generality of this method, and its straightforward extensions to three and higher degree of freedom systems, we have a 
unified framework to identify the dynamical mechanisms of transition in the 
presence of dissipation. 
While we  only considered two degree of freedom systems, in future work, higher dimensional systems will be considered. 

\section*{Acknowledgments}
This work was supported in part by the National Science Foundation under award 1537349.
The authors would like to thank Mingwu Li for the discussions on COCO and Hinke M Osinga for stimulating discussion on the early version of the draft. We also thank Harry Dankowicz and Jan Sieber for hosting ``Advanced Summer School on Continuation Methods for Nonlinear Problems'' at UIUC in 2018 from which the authors got to know COCO.

\noindent{\bf Appendices}

\begin{appendices}
\section{Quadratic equation for the transition ellipsoid}
\label{Appendix_quadratic}
The form of the transition ellipsoid in \eqref{ellipsoid_fourmula} that mediates the transition in the dissipative system for the snap-through buckling of a shallow arch can be rewritten by the following  form: 
	\begin{equation}
	a_{\bar p_2} \left( \bar p_2^0 \right)^2 + b_{\bar p_2} \bar p_2^0 + c_{\bar p_2}=0,
	\label{ellipsoid formula}
	\end{equation}
	where $a_{\bar p_2}$, $b_{\bar p_2}$, and $c_{\bar p_2}$ are given by
	\begin{equation*}
	\begin{aligned}
	& a_{\bar p_2}= \frac{s_2^2}{2 \omega_p(c_x + \omega_p^2)^2}, \hspace{0.2in} b_{\bar p_2}=\frac{\lambda s_2^2 (1+k_p) (c_x - \lambda^2) [\bar q_2^0 - \bar q_1^0 (c_x + \omega_p^2)]}{\omega_p (k_p - 1) (c_x + \omega_p^2)^2 (\lambda^2 + \omega_p^2)},\\
	& c_{\bar p_2}= c_p - \frac{\lambda^2 s_2^2 (1+k_p)^2 (c_x - c_y)[\bar q_2^0 - \bar q_1^0 (c_x + \omega_p^2)]^2}{2 \omega_p (k_p-1)^2 (c_x + \omega_p^2)^2(\lambda^2 + \omega_p^2)},\\
	&c_p = \left(\sum\limits_{i=1}^{4} c_p^{(i)}\right)/ \left[2 \omega_p \left(k_p-1 \right)^2 \left(\lambda^2 + \omega_p^2 \right)^2 \right]- h,\\
	& c_p^{(1)}=2 k_p s_1^2 \lambda \omega_p \left[\bar q_2- \bar q_1 \left(c_x + \omega_p^2 \right) \right]^2, \\
	& c_p^{(2)}= 8 k_p s_2^2 \lambda^2 \omega_p^2 \bar q_1 \left(c_x \bar q_1 - \bar q_2 \right),\\
	& c_p^{(3)}=s_2^2 \lambda^2 \left(1+k_p \right)^2 \left[\left(c_x \bar q_1- \bar q_2 \right)^2+ \bar q_1^2\omega_p^4 \right],\\
	& c_p^{(4)}=s_2^2 \omega_p^2 \left(k_p -1 \right)^2 \left[\left(c_x \bar q_1 - \bar q_2) \right)^2+ \bar q_1^2 \lambda^4 \right].
	\end{aligned}
	\end{equation*}

\end{appendices}

\bibliographystyle{shane-unsrt} 
\bibliography{Jun,Jun_general,ross_refs3}

\begin{thebibliography}{10}
\providecommand{\url}[1]{{\tt #1}}
\providecommand{\urlprefix}{URL }

\bibitem{zhong2018tube}
Zhong, J., Virgin, L.~N. and Ross, S.~D. [2018] A tube dynamics perspective
  governing stability transitions: An example based on snap-through buckling.
\newblock {\em International Journal of Mechanical Sciences\/} {\bf
  149}:413--428.

\bibitem{collins2012isomerization}
Collins, P., Ezra, G.~S. and Wiggins, S. [2012] Isomerization dynamics of a
  buckled nanobeam.
\newblock {\em Physical Review E\/} {\bf 86}(5):056218.

\bibitem{OzDeMeMa1990}
{Ozorio de Almeida}, A.~M., {De Leon}, N., Mehta, M.~A. and Marston, C.~C.
  [1990] Geometry and dynamics of stable and unstable cylinders in
  {H}amiltonian systems.
\newblock {\em Physica D\/} {\bf 46}:265--285.

\bibitem{DeMeTo1991}
{De Leon}, N., Mehta, M.~A. and Topper, R.~Q. [1991] Cylindrical manifolds in
  phase space as mediators of chemical reaction dynamics and kinetics. {I}.
  {T}heory.
\newblock {\em J. Chem. Phys.\/} {\bf 94}:8310--8328.

\bibitem{wiggins2001impenetrable}
Wiggins, S., Wiesenfeld, L., Jaff{\'e}, C. and Uzer, T. [2001] Impenetrable
  barriers in phase-space.
\newblock {\em Physical Review Letters\/} {\bf 86}(24):5478.

\bibitem{uzer2002geometry}
Uzer, T., Jaff{\'e}, C., Palaci{\'a}n, J., Yanguas, P. and Wiggins, S. [2002]
  The geometry of reaction dynamics.
\newblock {\em Nonlinearity\/} {\bf 15}(4):957.

\bibitem{gabern2005theory}
Gabern, F., Koon, W.~S., Marsden, J.~E. and Ross, S.~D. [2005] Theory and
  computation of non-RRKM lifetime distributions and rates in chemical systems
  with three or more degrees of freedom.
\newblock {\em Physica D: Nonlinear Phenomena\/} {\bf 211}(3-4):391--406.

\bibitem{jaffe2002statistical}
Jaff{\'e}, C., Ross, S.~D., Lo, M.~W., Marsden, J., Farrelly, D. and Uzer, T.
  [2002] Statistical theory of asteroid escape rates.
\newblock {\em Physical Review Letters\/} {\bf 89}(1):011101.

\bibitem{KoMaRoLoSc2004}
Koon, W.~S., Marsden, J.~E., Ross, S.~D., Lo, M.~W. and Scheeres, D.~J. [2004]
  Geometric mechanics and the dynamics of asteroid pairs.
\newblock {\em Annals of the New York Academy of Sciences\/} {\bf 1017}:11--38.

\bibitem{onozaki2017tube}
Onozaki, K., Yoshimura, H. and Ross, S.~D. [2017] Tube dynamics and low energy
  Earth--Moon transfers in the 4-body system.
\newblock {\em Advances in Space Research\/} {\bf 60}(10):2117--2132.

\bibitem{sequeira2018manifestation}
Sequeira, D., Wang, X.-S. and Mann, B. [2018] On the manifestation of
  coexisting nontrivial equilibria leading to potential well escapes in an
  inhomogeneous floating body.
\newblock {\em Physica D: Nonlinear Phenomena\/} {\bf 365}:80--90.

\bibitem{NaRo2017}
Naik, S. and Ross, S.~D. [2017] Geometry of escaping dynamics in nonlinear ship
  motion.
\newblock {\em Communications in Nonlinear Science and Numerical Simulation\/}
  {\bf 47}:48 -- 70.

\bibitem{zhong2020geometry}
Zhong, J. and Ross, S.~D. [2020] Geometry of escape and transition dynamics in
  the presence of dissipative and gyroscopic forces in two degree of freedom
  systems.
\newblock {\em Communications in Nonlinear Science and Numerical Simulation\/}
  {\bf 82}:105033.

\bibitem{meiss2007differential}
Meiss, J.~D. [2007] {\em Differential dynamical systems\/}, vol.~14.
\newblock Siam.

\bibitem{wiggins2003introduction}
Wiggins, S. [2003] {\em Introduction to applied nonlinear dynamical systems and
  chaos\/}, vol.~2.
\newblock Springer Science \& Business Media.

\bibitem{perko2013differential}
Perko, L. [2013] {\em Differential equations and dynamical systems\/}, vol.~7.
\newblock Springer Science \& Business Media.

\bibitem{Moser1958}
Moser, J. [1958] On the generalization of a theorem of {L}iapunov.
\newblock {\em Comm. Pure Appl. Math.\/} {\bf 11}:257--271.

\bibitem{Moser1973}
Moser, J. [1973] {\em Stable and Random Motions in Dynamical Systems with
  Special Emphasis on Celestial Mechanics\/}.
\newblock Princeton University Press.

\bibitem{krauskopf2006survey}
Krauskopf, B., Osinga, H.~M., Doedel, E.~J., Henderson, M.~E., Guckenheimer,
  J., Vladimirsky, A., Dellnitz, M. and Junge, O. [2006] A survey of methods
  for computing (un) stable manifolds of vector fields.
\newblock In {\em Modeling And Computations In Dynamical Systems: In
  Commemoration of the 100th Anniversary of the Birth of John von Neumann\/},
  67--95. World Scientific.

\bibitem{parker2012practical}
Parker, T.~S. and Chua, L. [2012] {\em Practical numerical algorithms for
  chaotic systems\/}.
\newblock Springer Science \& Business Media.

\bibitem{krauskopf2003computing}
Krauskopf, B. and Osinga, H.~M. [2003] Computing geodesic level sets on global
  (un) stable manifolds of vector fields.
\newblock {\em SIAM Journal on Applied Dynamical Systems\/} {\bf
  2}(4):546--569.

\bibitem{osinga2018understanding}
Osinga, H.~M. [2018] Understanding the geometry of dynamics: The stable
  manifold of the Lorenz system.
\newblock {\em Journal of the Royal Society of New Zealand\/} {\bf
  48}(2-3):203--214.

\bibitem{dellnitz1997subdivision}
Dellnitz, M. and Hohmann, A. [1997] A subdivision algorithm for the computation
  of unstable manifolds and global attractors.
\newblock {\em Numerische Mathematik\/} {\bf 75}(3):293--317.

\bibitem{dellnitz1996computation}
Dellnitz, M. and Hohmann, A. [1996] The computation of unstable manifolds using
  subdivision and continuation.
\newblock In {\em Nonlinear dynamical systems and chaos\/}, 449--459. Springer.

\bibitem{madrid2009distinguished}
Madrid, J.~J. and Mancho, A.~M. [2009] Distinguished trajectories in time
  dependent vector fields.
\newblock {\em Chaos: An Interdisciplinary Journal of Nonlinear Science\/} {\bf
  19}(1):013111.

\bibitem{mendoza2010hidden}
Mendoza, C. and Mancho, A.~M. [2010] Hidden geometry of ocean flows.
\newblock {\em Physical review letters\/} {\bf 105}(3):038501.

\bibitem{naik2019finding}
Naik, S. and Wiggins, S. [2019] Finding normally hyperbolic invariant manifolds
  in two and three degrees of freedom with H{\'e}non-Heiles-type potential.
\newblock {\em Physical Review E\/} {\bf 100}(2):022204.

\bibitem{mancho2013lagrangian}
Mancho, A.~M., Wiggins, S., Curbelo, J. and Mendoza, C. [2013] Lagrangian
  descriptors: A method for revealing phase space structures of general time
  dependent dynamical systems.
\newblock {\em Communications in Nonlinear Science and Numerical Simulation\/}
  {\bf 18}(12):3530--3557.

\bibitem{dankowicz2013recipes}
Dankowicz, H. and Schilder, F. [2013] {\em Recipes for continuation\/},
  vol.~11.
\newblock SIAM.

\bibitem{zhong2016analysis}
Zhong, J., Fu, Y., Chen, Y. and Li, Y. [2016] Analysis of nonlinear dynamic
  responses for functionally graded beams resting on tensionless elastic
  foundation under thermal shock.
\newblock {\em Composite Structures\/} {\bf 142}:272--277.

\bibitem{WiVi2016}
Wiebe, R. and Virgin, L.~N. [2016] On the experimental identification of
  unstable static equilibria.
\newblock {\em Proceedings of the Royal Society of London A: Mathematical,
  Physical and Engineering Sciences\/} {\bf 472}(2190):20160172.

\bibitem{zhong2020differential}
Zhong, J. and Ross, S.~D. [2020] Differential correction and arc-length
  continuation applied to boundary value problems: examples based on
  snap-through of circular arches and spherical shells.
\newblock {\em viXra\/} .

\bibitem{virgin2017geometric}
Virgin, L., Guan, Y. and Plaut, R. [2017] On the geometric conditions for
  multiple stable equilibria in clamped arches.
\newblock {\em International Journal of Non-Linear Mechanics\/} {\bf 92}:8--14.

\bibitem{harvey2015coexisting}
Harvey~Jr, P. and Virgin, L. [2015] Coexisting equilibria and stability of a
  shallow arch: Unilateral displacement-control experiments and theory.
\newblock {\em International Journal of Solids and Structures\/} {\bf
  54}:1--11.

\bibitem{murrell1968symmetries}
Murrell, J.~N. and Laidler, K.~J. [1968] Symmetries of activated complexes.
\newblock {\em Transactions of the Faraday Society\/} {\bf 64}:371--377.

\bibitem{DeLi1994}
{De Leon}, N. and Ling, S. [1994] Simplification of the transition state
  concept in reactive island theory: {A}pplication to the
  {HCN$\rightleftharpoons$CNH} isomerization.
\newblock {\em J. Chem. Phys.\/} {\bf 101}:4790--4802.

\bibitem{delavega1982role}
{De la {V}ega}, J.~R. [1982] Role of symmetry in the tunneling of the proton in
  double-minimum potentials.
\newblock {\em Accounts of Chemical Research\/} {\bf 15}(6):185--191.

\bibitem{minyaev1994reaction}
Minyaev, R.~M. [1994] Reaction path as a gradient line on a potential energy
  surface.
\newblock {\em International Journal of Quantum Chemistry\/} {\bf
  49}(2):105--127.

\bibitem{smedarchina2007correlated}
Smedarchina, Z., Siebrand, W. and Fern{\'a}ndez-Ramos, A. [2007] Correlated
  double-proton transfer. {I}. {T}heory.
\newblock {\em Journal of Chemical Physics\/} {\bf 127}(17):174513.

\bibitem{accardi2010synchronous}
Accardi, A., Barth, I., Kühn, O. and Manz, J. [2010] From synchronous to
  sequential double proton transfer: Quantum dynamics simulations for the model
  Porphine.
\newblock {\em The Journal of Physical Chemistry A\/} {\bf
  114}(42):11252--11262.

\bibitem{ezra2009phase}
Ezra, G.~S. and Wiggins, S. [2009] Phase-space geometry and reaction dynamics
  near index 2 saddles.
\newblock {\em Journal of Physics A: Mathematical and Theoretical\/} {\bf
  42}(20):205101.

\bibitem{Greenwood2003}
Greenwood, D.~T. [2003] {\em Advanced Dynamics\/}.
\newblock Cambridge University Press.

\bibitem{Wiggins1994}
Wiggins, S. [1994] {\em Normally Hyperbolic Invariant Manifolds in Dynamical
  Systems\/}.
\newblock Springer-Verlag, New York.

\bibitem{McGehee1969}
McGehee, R. [1969] {\em Some homoclinic orbits for the restricted three-body
  problem\/}.
\newblock Ph.D. thesis, University of Wisconsin, Madison.

\bibitem{Conley1968}
Conley, C.~C. [1968] Low energy transit orbits in the restricted three-body
  problem.
\newblock {\em SIAM J. Appl. Math.\/} {\bf 16}:732--746.

\bibitem{seydel2009practical}
Seydel, R. [2009] {\em Practical bifurcation and stability analysis\/}, vol.~5.
\newblock Springer Science \& Business Media.

\bibitem{nayfeh2008nonlinear}
Nayfeh, A.~H. and Mook, D.~T. [2008] {\em Nonlinear oscillations\/}.
\newblock John Wiley \& Sons.

\bibitem{nayfeh2008applied}
Nayfeh, A.~H. and Balachandran, B. [2008] {\em Applied nonlinear dynamics:
  analytical, computational, and experimental methods\/}.
\newblock John Wiley \& Sons.

\bibitem{lau1981amplitude}
Lau, S.~L. and Cheung, Y.~K. [1981] Amplitude incremental variational principle
  for nonlinear vibration of elastic systems.
\newblock {\em Journal of Applied Mechanics\/} {\bf 48}(4):959--964.

\bibitem{fu2006analysis}
Fu, Y., Hong, J. and Wang, X. [2006] Analysis of nonlinear vibration for
  embedded carbon nanotubes.
\newblock {\em Journal of Sound and Vibration\/} {\bf 296}(4-5):746--756.

\bibitem{KoLoMaRo2011}
Koon, W.~S., Lo, M.~W., Marsden, J.~E. and Ross, S.~D. [2011] {\em Dynamical
  Systems, the Three-Body Problem and Space Mission Design\/}.
\newblock Marsden Books, ISBN 978-0-615-24095-4.

\bibitem{sundararajan1997dynamics}
Sundararajan, P. and Noah, S.~T. [1997] Dynamics of forced nonlinear systems
  using shooting/arc-length continuation method---application to rotor systems.
\newblock {\em Journal of Vibration and Acoustics\/} {\bf 119}(1):9--20.

\bibitem{krauskopf2007numerical}
Krauskopf, B., Osinga, H.~M. and Gal{\'a}n-Vioque, J. [2007] {\em Numerical
  continuation methods for dynamical systems\/}.
\newblock Springer.

\end{thebibliography}

\end{document}